\documentclass[12pt,a4paper]{article}
\usepackage[14pt]{extsizes}
\usepackage{iftex}
\ifPDFTeX
  \usepackage[utf8]{inputenc}
\fi
\usepackage[english]{babel}
\usepackage{amsmath,amsfonts,amssymb}
\usepackage{geometry}
\usepackage{bm}
\usepackage{booktabs}
\usepackage{float}
\usepackage{tikz}

\tikzset{
  satake/.style={line width=0.7pt, x=1.25cm, y=1.25cm, baseline={([yshift=-0.5ex]current bounding box.center)}},
  ev/.style={circle, draw, fill=white, inner sep=0pt, minimum size=8pt},
  od/.style={rectangle, draw, fill=white, inner sep=0pt, minimum size=7.4pt},
  lab/.style={font=\scriptsize, inner sep=1.5pt},
  dbl/.style={double, double distance=2.2pt},
  tauarc/.style={line width=0.6pt},
}
\newcommand{\SatSOthree}{\tikz[satake]{
  \node[ev, label={[lab]below:$\alpha_1$}] at (0,0) {};}}
\newcommand{\SatSOfive}{\tikz[satake]{
  \draw[dbl] (0,0) -- (1.1,0);
  \draw (0.49,0.13) -- (0.61,0) -- (0.49,-0.13);
  \node[ev, label={[lab]below:$\alpha_1$}] at (0,0) {};
  \node[ev, label={[lab]below:$\alpha_2$}] at (1.1,0) {};}}
\newcommand{\SatOSPonetwo}{\tikz[satake]{
  \node[od, label={[lab]below:$\alpha_1$}] at (0,0) {};}}
\newcommand{\SatOSPonefour}{\tikz[satake]{
  \draw[dbl] (0,0) -- (1.1,0);
  \draw (0.49,0.13) -- (0.61,0) -- (0.49,-0.13);
  \node[ev, label={[lab]below:$\alpha_1$}] at (0,0) {};
  \node[od, label={[lab]below:$\alpha_2$}] at (1.1,0) {};}}
\newcommand{\SatSOfour}{\tikz[satake]{
  \node[ev, label={[lab]above:$\alpha_1$}] at (0,0.6) {};
  \node[ev, label={[lab]below:$\alpha_2$}] at (0,-0.6) {};}}
\newcommand{\SatSOsix}{\tikz[satake]{
  \draw (0,0) -- (1,0.6) (0,0) -- (1,-0.6);
  \node[ev, label={[lab]left:$\alpha_1$}] at (0,0) {};
  \node[ev, label={[lab]above:$\alpha_2$}] (b) at (1,0.6) {};
  \node[ev, label={[lab]below:$\alpha_3$}] (c) at (1,-0.6) {};
  \draw[tauarc] (b.east) to[out=-40, in=40, looseness=1.1] (c.east);}}
\newcommand{\SatOSPtwotwo}{\tikz[satake]{
  \draw[dbl] (0,0.6) -- (0,-0.6);
  \node[od, label={[lab]above:$\alpha_1$}] (b) at (0,0.6) {};
  \node[od, label={[lab]below:$\alpha_2$}] (c) at (0,-0.6) {};
  \draw[tauarc] (b.east) to[out=-40, in=40, looseness=1.1] (c.east);}}
\newcommand{\SatOSPtwofour}{\tikz[satake]{
  \draw (0,0) -- (1,0.6) (0,0) -- (1,-0.6);
  \draw[dbl] (1,0.6) -- (1,-0.6);
  \node[ev, label={[lab]left:$\alpha_1$}] at (0,0) {};
  \node[od, label={[lab]above:$\alpha_2$}] (b) at (1,0.6) {};
  \node[od, label={[lab]below:$\alpha_3$}] (c) at (1,-0.6) {};
  \draw[tauarc] (b.east) to[out=-40, in=40, looseness=1.1] (c.east);}}

\allowdisplaybreaks
\newtheorem{theorem}{Theorem}
\newtheorem{propn}[theorem]{Proposition}
\newcommand{\g}{\mathfrak{g}}
\newcommand{\kk}{\mathfrak{k}}
\newcommand{\so}{\mathfrak{so}}
\newcommand{\osp}{\mathfrak{osp}}
\newcommand{\Uq}{U_q}
\newcommand{\ta}{\tilde\alpha}
\newcommand{\ot}{\otimes}
\newcommand{\diag}{\operatorname{diag}}
\newcommand{\spn}{\operatorname{span}}
\newcommand{\hd}[2][c]{\tabular[b]{@{}#1@{}}#2\endtabular}

\title{Orthosymplectic $K$-matrices in small dimensions}
\author{E.~Ryzhkin${}^{1}$, V.~Stukopin${}^{1,2}$
\vspace{10pt}\\
\small ${}^{1}$ Moscow Institute of Physics and Technology\\
\small
${}^{2}$ Southern Mathematical Institute\\
\small e-mail:
ryzhkin.ea@phystech.edu, stukopin.va@mipt.ru\\
}

\begin{document}
\sloppy
\maketitle

\begin{abstract}
Using elementary methods, we describe coideal subalgebras of quantum orthosymplectic algebras of small dimensions, classify the corresponding $K$-matrices and find explicitly the centralizers of these subalgebras in the vector representation. The main attention is paid to generators with a Cartan term. For them, at generic values of the parameters, we find the centralizers and all $K$-matrices and describe the decomposition of the vector representation over the coideal subalgebra. We find all special values of the parameters of two types and, at a generic point of each set of special values, describe the decomposition of the vector representation and all $K$-matrices. At special values of the first type the dimension of the centralizer grows. One hypersurface of such values is related to Cartan terms, and on it one of the generators is given by the same formula as in the quantum electrical algebra. At special values of the second type the dimension of the centralizer is the same, but nonzero nilpotent elements appear in it. The $K$-matrices on the sets of special values of the second type are obtained from those already found by substituting the condition, and there are no others. We describe the classical limit, including the limit of special values. We formulate an observation relating the admissibility of a Cartan term to the first homology of the classical subalgebra.
\end{abstract}

\section{Introduction}

The reflection equation (RE) plays a fundamental role in the theory of quantum groups and quantum integrable systems. It is as important an ingredient of the general theory as the quantum Yang--Baxter equation. The reflection equation arises in connection with quantum integrable systems with boundary conditions. The reflection equation algebras arising in this way turn out to be closely related to the theory of coideal subalgebras, initiated by the fundamental work of Letzter~\cite{L} and developed over a long time (see, for example, \cite{K}, \cite{K1}, which reflect the progress achieved during the first decade of development of this field). The analogous theory for quantum superalgebras started to develop relatively recently and is now an area of active research (see~\cite{AMS}, \cite{SW}). It should be noted that there are several versions of the reflection equation: for example, along with the ordinary version there is also a twisted one. In recent works, the solutions of the reflection equation~\cite{AMS1} and of the twisted reflection equation mentioned above~\cite{AM} were classified, in both cases related to the quantum superalgebra $\mathfrak{gl}(m|n)$. The case of orthosymplectic superalgebras, especially in small dimensions, is rather specific, and its study is of current interest because, on the one hand, effects that disappear as the dimension grows often occur here and, on the other hand, the case of small dimensions is, as a rule, important for applications in the theory of quantum integrable systems. It should be said that the choice of different $K$-matrices in the theory of integrable systems often corresponds to the choice of different boundary conditions. In this paper we describe the solutions of the reflection equation for orthosymplectic superalgebras of small dimensions. These solutions do not always fit into the general classification that holds for sufficiently large dimensions. The question of such a general classification for orthosymplectic Lie superalgebras is considered in~\cite{AMS}, but there series of types of coideal subalgebras and of the corresponding $K$-matrices are described, and there is no detailed analysis of the solutions for different values of the parameters. In this paper we give a detailed description of the solutions of the reflection equation and of the coideal subalgebras. We obtain this description using only elementary methods. Along the way we note the appearance of interesting examples of coideal subalgebras at special values of the parameters. Thus we consider examples of orthosymplectic analogues of electrical Lie algebras (orthosymplectic Lie superalgebras) (see \cite{AMS2}, \cite{LP}, \cite{BGG}) in the case of small dimensions, which arise as classical limits of coideal subalgebras at special values of the parameters of the generators of the coideal subalgebras. Note that, besides using elementary methods, we try to explain the constructions used and the results obtained in sufficient detail. We hope that this work will be useful for further study of the structure of coideal subalgebras in the general case, that the results obtained here in particular cases for exceptional values of the parameters will help in further research, and that they will also be useful in the theory of quantum integrable systems with boundary conditions.

We consider the quantized universal enveloping algebras $\Uq(\g)$ of the Lie algebras $\so(3)$, $\so(5)$, $\so(4)$, $\so(6)$ and of the Lie superalgebras $\osp(1|2)$, $\osp(1|4)$,
$\osp(2|2)$, $\osp(2|4)$ in the vector representation $V$. A Satake diagram defines a subalgebra $\kk\subset\g$. A matrix $K\in\mathrm{End}(V)$ is called \emph{invariant} if it commutes with the representation of the coideal subalgebra $\mathcal B=\Uq(\kk)\subset\Uq(\g)$. In this paper a \emph{$K$-matrix} is
an even invertible invariant matrix satisfying the reflection equation.

The generators of the coideal subalgebra for a white root $\alpha$ have the form $X_\alpha=T_\alpha e_\alpha+c_\alpha F_{\ta}$
(Section~\ref{ss:gen}). One can add a Cartan term to them:
\[
X^{(s)}_\alpha=X_\alpha+s\,(q^{h_{\ta}}-1),\qquad s\in\mathbb C,
\]
which corresponds to the parameters $s$ of quantum symmetric pairs~\cite{K1}. The subalgebra then remains a left coideal
(Section~\ref{ss:Xs}). The case $s\neq0$ is the main one in this paper.

The following results are obtained.
\begin{enumerate}
\item We describe the representations of the coideal subalgebras $\mathcal B$ of the algebras mentioned above and all the corresponding $K$-matrices
at generic values of the parameters (Theorem~\ref{th:main}).
\item For generators with a Cartan term, at generic values of the parameters, we find the centralizers and all $K$-matrices (Theorem~\ref{th:cartan}).
Nonscalar solutions exist for $\so(3)$, $\so(4)$, for the node $\alpha_1$ in $\so(5)$ and $\osp(1|4)$ and, for
$c_2=c_3$, for the node $\alpha_1$ in $\so(6)$ and $\osp(2|4)$. In the remaining cases the centralizer consists of scalar
matrices.
\item Cartan terms are considered on all sets of nodes. Apart from the condition $c_2=c_3$ for $\so(6)$ and $\osp(2|4)$,
there is exactly one hypersurface of special values of the first type on which the centralizer is larger than in the general case:
$\so(5)$ with terms on both nodes and $c_1=-(q-1)^2s^2$. On it the generator $X^{(s)}_1$ has the form of a generator of the quantum electrical algebra~\cite{AMS2}, the $K$-matrix
has three eigenvalues of multiplicities $2$, $2$, $1$ (rather than $3$, $1$, $1$), and the classical subalgebra is isomorphic to
$\mathfrak{gl}(2)$, as for the term on $\alpha_1$ only, but $V$ decomposes over it as $2+2+1$ (rather than $3+1+1$).
\item In the cases of Table~\ref{t:dec} the centralizer at a generic point is commutative and semisimple, and $V$ decomposes
over $\mathcal B$ without multiplicities. We find all special values of the second type, at which the dimension of the centralizer is the same but nonzero nilpotent elements
appear in it (Table~\ref{t:typeII}). To find them we use the trace form on the
centralizer. Inside these hypersurfaces there are special values only for $\so(4)$, at the intersection of two
hypersurfaces, where $V$ is indecomposable. On the sets of special values of the second type the $K$-matrices are obtained from
those already found by substituting the condition, and there are no others (Proposition~\ref{p:typeII}).
\item We describe the classical limit: the subalgebras $\kk_\sigma\subset\g$, their derived algebras and first homologies, and also the
limit of special values (Table~\ref{t:cl3}).
We formulate an observation on the condition under which a Cartan term admits a nonscalar $K$-matrix.
\end{enumerate}
For $q\neq1$ the $K$-matrices were found in two steps. First we found the general form of an even matrix commuting with the generators
of the coideal subalgebra, i.e., the centralizer in the representation, and then the solutions of the reflection equation were singled out among its matrices.


\section{Preliminaries}
The definitions and notation of this paper follow~\cite{AMS}. For general facts on quantum groups see~\cite{KS}, and general facts on superalgebras can be found in the book~\cite{CW}. For the reader's convenience we recall the main definitions and facts used below.

\subsection{Quantum superalgebra}
The algebra $\Uq(\g)$ is generated by the elements $e_\alpha$, $f_\alpha$, $q^{\pm h_\alpha}$, where $\alpha$ runs over the set $\Pi$ of simple roots. By $|\Pi|$ we denote the number of simple roots.
Parity: $|q^{h_\alpha}|=0$, $|e_\alpha|=|f_\alpha|=|\alpha|$. These generators satisfy, in particular, the following relations, which involve the entries of the symmetrized Cartan matrix $C=\bigl((\alpha_i,\alpha_j)\bigr)_{i,j=1}^{|\Pi|}$, the scalar products of the simple roots:
\begin{equation}
\begin{gathered}
q^{h_\alpha}e_\beta q^{-h_\alpha}=q^{(\alpha,\beta)}e_\beta,\qquad
q^{h_\alpha}f_\beta q^{-h_\alpha}=q^{-(\alpha,\beta)}f_\beta,\\
e_\alpha f_\beta-(-1)^{|\alpha||\beta|}f_\beta e_\alpha=0\ \ (\alpha\neq\beta),
\end{gathered}
\end{equation}
the coproduct is given by the formulas
\begin{equation}\label{eq:cop}
\begin{gathered}
\Delta(e_\alpha)=q^{h_\alpha}\ot e_\alpha+e_\alpha\ot1,\qquad
\Delta(f_\alpha)=f_\alpha\ot q^{-h_\alpha}+1\ot f_\alpha,\\
\Delta(q^{h_\alpha})=q^{h_\alpha}\ot q^{h_\alpha},
\end{gathered}
\end{equation}
and the graded tensor product obeys the sign rule
\begin{equation}\label{eq:koszul}
(a\ot b)(c\ot d)=(-1)^{|b||c|}\,ac\ot bd .
\end{equation}
Note that the coproduct is a homomorphism of superalgebras: $\Delta(ab)=\Delta(a)\Delta(b)$. Moreover, we have omitted here the additional Serre relations (see \cite{Y}), which do not play an important role in what follows.

\subsection{Vector representation}\label{ss:rep}
The basis $v_1,\dots,v_N$ of the space $V$ consists of weight vectors. Their weights are $\pm\varepsilon_i$, $\pm\delta_j$ and, for odd $N$, the zero weight. The vector $v_i$ has parity $|i|\in\{0,1\}$: the vectors of weights $\pm\delta_j$ are odd, the others are even.
Set $i'=N+1-i$. The elementary matrices $E_{ij}$ (whose only nonzero entry, equal to 1, is at the intersection of the $i$-th row and the $j$-th column) have parity $|i|+|j|$. The element $q^{h_\alpha}$ acts on a vector of weight $\nu$
by multiplication by $q^{(\alpha,\nu)}$, where $(\varepsilon_i,\varepsilon_j)=\delta_{ij}$, $(\delta_i,\delta_j)=-\delta_{ij}$, $(\varepsilon_i,\delta_j)=0$.
The gradings are symmetric: $|i'|=|i|$. The data of the representations are taken from~\cite[Sections~2.2, 2.3]{AMS}. The parities of the basis vectors and the simple roots are collected in
Table~\ref{t:data}.

\begin{table}[H]
\centering\small
\begin{tabular}{@{}lll@{}}
\toprule
$\g$ & parities & simple roots\\
\midrule
$\so(3)$ & $(0,0,0)$ & $\varepsilon_1$\\
$\so(5)$ & $(0,0,0,0,0)$ & $\varepsilon_1-\varepsilon_2$,\ $\varepsilon_2$\\
$\so(4)$ & $(0,0,0,0)$ & $\varepsilon_1-\varepsilon_2$,\ $\varepsilon_1+\varepsilon_2$\\
$\so(6)$ & $(0,0,0,0,0,0)$ & $\varepsilon_1-\varepsilon_2$,\ $\varepsilon_2-\varepsilon_3$,\ $\varepsilon_2+\varepsilon_3$\\
$\osp(1|2)$ & $(1,0,1)$ & $\delta_1$ (odd)\\
$\osp(1|4)$ & $(1,1,0,1,1)$ & $\delta_1-\delta_2$,\ $\delta_2$ (odd)\\
$\osp(2|2)$ & $(1,0,0,1)$ & $\delta_1-\varepsilon_1$,\ $\delta_1+\varepsilon_1$ (odd)\\
$\osp(2|4)$ & $(1,1,0,0,1,1)$ & $\delta_1-\delta_2$,\ $\delta_2-\varepsilon_1$ (odd),\ $\delta_2+\varepsilon_1$ (odd)\\
\bottomrule
\end{tabular}
\caption{Parities of the basis vectors and simple roots.}\label{t:data}
\end{table}

Matrices of the generators (for $\so(N)$ everywhere $f_i=e_i^{t}$):
\begin{itemize}
\item $\so(3)$: $e_1=E_{12}-E_{23}$, $q^{h_1}=\diag(q,1,q^{-1})$,
\item $\so(5)$: $e_1=E_{12}-E_{45}$, $e_2=E_{23}-E_{34}$, $q^{h_1}=\diag(q,q^{-1},1,q,q^{-1})$, $q^{h_2}=\diag(1,q,1,q^{-1},1)$,
\item $\so(4)$: $e_1=E_{12}-E_{34}$, $e_2=E_{13}-E_{24}$, $q^{h_1}=\diag(q,q^{-1},q,q^{-1})$, $q^{h_2}=\diag(q,q,q^{-1},q^{-1})$,
\item $\so(6)$: $e_1=E_{12}-E_{56}$, $e_2=E_{23}-E_{45}$, $e_3=E_{24}-E_{35}$,
$q^{h_1}=\diag(q,q^{-1},1,1,q,q^{-1})$, $q^{h_2}=\diag(1,q,q^{-1},q,q^{-1},1)$, $q^{h_3}=\diag(1,q,q,q^{-1},q^{-1},1)$,
\item $\osp(1|2)$: $e_1=q^{-1}E_{12}+E_{23}$, $f_1=-qE_{21}+E_{32}$, $q^{h_1}=\diag(q^{-1},1,q)$,
\item $\osp(1|4)$: $e_1=E_{12}-E_{45}$, $f_1=E_{21}-E_{54}$, $q^{h_1}=\diag(q^{-1},q,1,q^{-1},q)$,
$e_2=q^{-1}E_{23}+E_{34}$, $f_2=-qE_{32}+E_{43}$, $q^{h_2}=\diag(1,q^{-1},1,q,1)$,
\item $\osp(2|2)$: $e_1=q^{-1}E_{12}+E_{34}$, $f_1=-qE_{21}+E_{43}$, $q^{h_1}=\diag(q^{-1},q^{-1},q,q)$,
$e_2=q^{-1}E_{13}+E_{24}$, $f_2=-qE_{31}+E_{42}$, $q^{h_2}=\diag(q^{-1},q,q^{-1},q)$,
\item $\osp(2|4)$: $e_1=E_{12}-E_{56}$, $f_1=E_{21}-E_{65}$, $q^{h_1}=\diag(q^{-1},q,1,1,q^{-1},q)$,
$e_2=q^{-1}E_{23}+E_{45}$, $f_2=-qE_{32}+E_{54}$, $q^{h_2}=\diag(1,q^{-1},q^{-1},q,q,1)$,
$e_3=q^{-1}E_{24}+E_{35}$, $f_3=-qE_{42}+E_{53}$, $q^{h_3}=\diag(1,q^{-1},q,q^{-1},q,1)$.
\end{itemize}
In these representations the following relation holds:
\[
e_\alpha f_\alpha-(-1)^{|\alpha|}f_\alpha e_\alpha=\pm\frac{q^{h_\alpha}-q^{-h_\alpha}}{q-q^{-1}} .
\]
On the left there is a commutator if the root $\alpha$ is even and an anticommutator if it is odd.

\subsection{Satake diagram and the subalgebra $\kk$}\label{ss:satake}
Let us recall the main definitions related to the theory of Lie superalgebras and their subalgebras. Basic Lie superalgebras, like simple Lie algebras, are described by Dynkin diagrams whose vertices correspond to simple roots (see~\cite{CW}). We draw even simple roots as circles and odd ones as squares.
A Satake diagram associated with a basic Lie superalgebra $\g$ is the Dynkin diagram of $\g$ in which some of the nodes are painted black and some
pairs of white nodes are connected by an arc (which colorings and arcs are admissible is described in~\cite[Section~4]{AMS}). The arcs define a permutation $\tau$ of the simple roots: two roots
connected by an arc are interchanged, the others stay in place. The diagram defines a subalgebra $\kk\subset\g$ by the following rules.
\begin{itemize}
\item If a root $\alpha$ is black, then $e_\alpha$, $f_\alpha$ and $h_\alpha$ lie in $\kk$.
\item If a root $\alpha$ is white, then instead of $e_\alpha$ and $f_\alpha$ one takes only the combination of the raising and lowering elements
\begin{equation}\label{eq:cl}
x_\alpha=e_\alpha+c_\alpha f_{\ta}
\end{equation}
and the Cartan element $h_{\ta}-h_\alpha$. Here $c_\alpha\neq0$ is a number called the mixing parameter,
and $\ta$ is the positive root assigned to the root $\alpha$.
\end{itemize}
The subalgebra $\kk$ is generated by all these elements.

The root $\ta$ is constructed as follows. If there are no black nodes, then $\ta=\tau(\alpha)$: this is the root $\alpha$ itself or the other end of the arc. If there are black nodes, then $\ta=w(\tau(\alpha))$, where $w$ is
the Weyl operator of the black part of the diagram~\cite[Section~4.1]{AMS}, an analogue of the longest element of the Weyl group. The root $\ta$
may turn out to be not simple but a sum of several simple roots (a composite root). In the example of Section~\ref{ss:comp} there is one black
root $\beta$, it is even, and $w$ is the reflection
\[
s_\beta(\gamma)=\gamma-\frac{2(\gamma,\beta)}{(\beta,\beta)}\,\beta .
\]

Examples. For $\so(3)$ the diagram consists of one white node, and $\kk$ is generated by one element
$x_1=e_1+c_1f_1$. For $\so(6)$ with a twisted tail the arc connects the nodes $\alpha_2$ and $\alpha_3$, therefore
$x_2=e_2+c_2f_3$, $x_3=e_3+c_3f_2$, and $\kk$ contains $h_3-h_2$.

\begin{figure}[H]
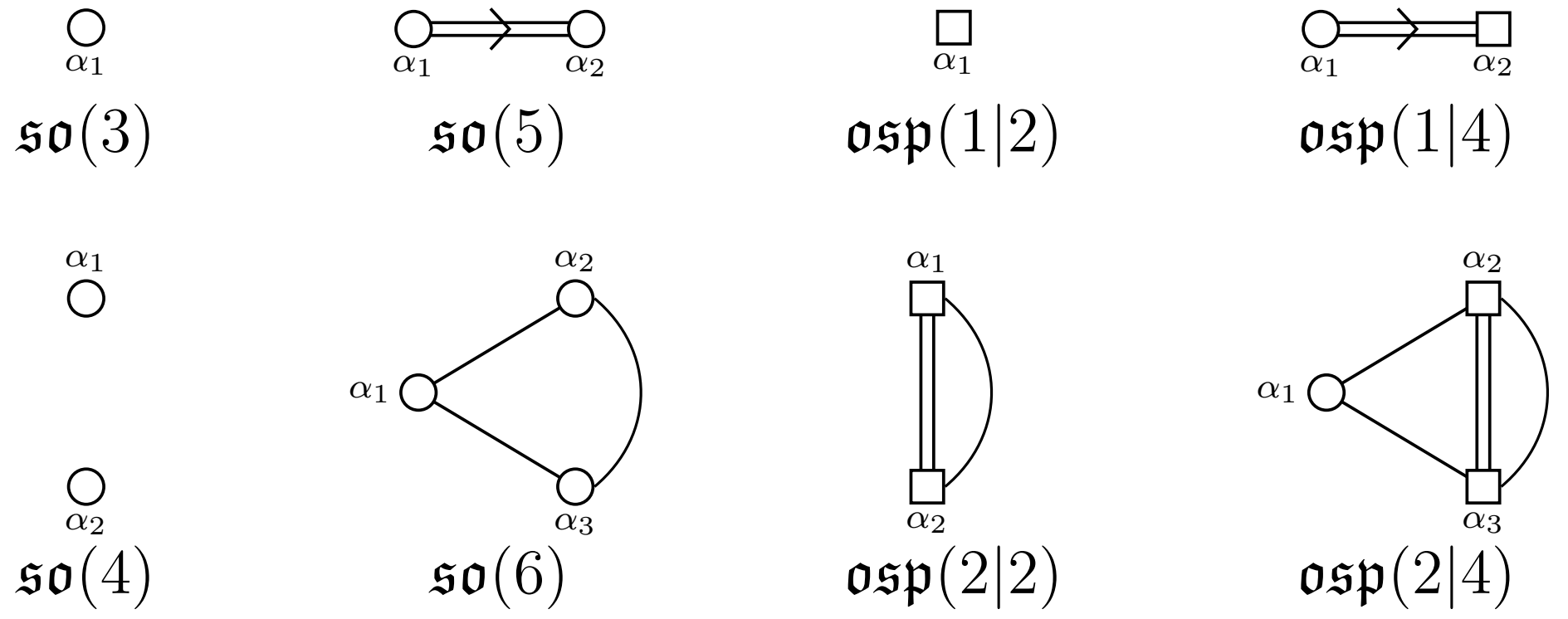

\centering
\begin{tabular}{@{}c@{\hspace{3em}}c@{\hspace{3em}}c@{\hspace{3em}}c@{}}
\SatSOthree & \SatSOfive & \SatOSPonetwo & \SatOSPonefour\\[4pt]
$\so(3)$ & $\so(5)$ & $\osp(1|2)$ & $\osp(1|4)$\\[16pt]
\SatSOfour & \SatSOsix & \SatOSPtwotwo & \SatOSPtwofour\\[4pt]
$\so(4)$ & $\so(6)$ & $\osp(2|2)$ & $\osp(2|4)$
\end{tabular}
\caption{Satake diagrams of the algebras under consideration. A circle is an even simple root, a square is an odd one.
An arc connects the roots interchanged by $\tau$.}\label{f:satake}
\end{figure}

In all eight cases under consideration all nodes are white (Figure~\ref{f:satake}). The permutation $\tau$ is trivial for $\so(3)$, $\so(5)$,
$\so(4)$, $\osp(1|2)$, $\osp(1|4)$ and interchanges the two roots of the tail for $\so(6)$, $\osp(2|2)$, $\osp(2|4)$
(a ``twisted tail'').

\subsection{Left coideal, $R$-matrix, reflection equation}\label{ss:RE}
A subalgebra $\mathcal B\subset\Uq(\g)$ is a \emph{left coideal} if $\Delta(\mathcal B)\subset\Uq(\g)\ot\mathcal B$.
Since $\Delta$ is a homomorphism, it suffices to check this on the generators of $\mathcal B$.

The $R$-matrix (formula (2.2) in~\cite{AMS}), $\omega=q-q^{-1}$:
\begin{multline}\label{eq:R}
R=\sum_{i,j}q^{(-1)^{|j|}(\delta_{ij}-\delta_{ij'})}E_{ii}\ot E_{jj}\\
+\omega\sum_{i>j}\Bigl((-1)^{|j|}E_{ij}\ot E_{ji}-(-1)^{|i|+|j|+|i||j|}\kappa_i\kappa_jq^{\rho_i-\rho_j}E_{ij}\ot E_{i'j'}\Bigr).
\end{multline}
The numbers $\kappa_i$ and $\rho_i$ are given in Table~\ref{t:kr}.

\begin{table}[H]
\centering\small
\renewcommand{\arraystretch}{1.2}
\begin{tabular}{@{}lll@{}}
\toprule
$\g$ & $\kappa$ & $\rho$\\
\midrule
$\so(3)$ & $(1,1,1)$ & $(\tfrac12,\tfrac12,-\tfrac12)$\\
$\so(5)$ & all $1$ & $(\tfrac32,\tfrac12,\tfrac12,-\tfrac12,-\tfrac32)$\\
$\so(4)$ & all $1$ & $(1,0,0,-1)$\\
$\so(6)$ & all $1$ & $(2,1,0,0,-1,-2)$\\
$\osp(1|2)$ & $(-1,1,1)$ & $(-\tfrac12,\tfrac12,\tfrac12)$\\
$\osp(1|4)$ & $(-1,-1,1,1,1)$ & $(-\tfrac32,-\tfrac12,\tfrac12,\tfrac12,\tfrac32)$\\
$\osp(2|2)$ & $(-1,1,1,1)$ & $(0,0,0,0)$\\
$\osp(2|4)$ & $(-1,-1,1,1,1,1)$ & $(-1,0,0,0,0,1)$\\
\bottomrule
\end{tabular}
\caption{The numbers $\kappa_i$ and $\rho_i$ in formula~\eqref{eq:R}.}\label{t:kr}
\end{table}

Note that operators of the form $A\ot B$ act on $V\ot V$ by the rule
\[
(A\ot B)(v\ot w)=(-1)^{|B||v|}\,Av\ot Bw ,
\]
consistent with~\eqref{eq:koszul}. The graded permutation $P$ interchanges the factors in $V\ot V$:
\[
P(v_i\ot v_j)=(-1)^{|i||j|}\,v_j\ot v_i .
\]
The sign ``$-$'' occurs only when both vectors are odd. In terms of the elementary matrices $E_{ij}$ the operator is $P=\sum_{i,j}(-1)^{|j|}E_{ij}\ot E_{ji}$.
Indeed, $(E_{ij}\ot E_{ji})(v_j\ot v_i)=(-1)^{(|i|+|j|)|j|}\,v_i\ot v_j$, and together with the factor $(-1)^{|j|}$
this gives the sign $(-1)^{|i||j|}$. Using $P$, one defines $R_{21}=PR_{12}P$ (that is, the matrix $R$ with the tensor factors interchanged) and $S=PR$.
The resulting braid relation $S_{12}S_{23}S_{12}=S_{23}S_{12}S_{23}$ is equivalent to the Yang--Baxter equation for $R$.

The reflection equation $SK_2SK_2=K_2SK_2S$ (formula (3.3) in~\cite{AMS}) for an even matrix $K$, after the substitution $S=PR$,
takes the form
\begin{equation}\label{eq:RE}
R_{21}K_1R_{12}K_2=K_2R_{21}K_1R_{12},\qquad K_1=K\ot1,\ \ K_2=1\ot K .
\end{equation}

For a twisted tail, \cite[(3.5)]{AMS} considers the reflection equation twisted by the automorphism $\tau$:
\begin{equation}\label{eq:tRE}
R_{21}K_1R^{\tau_1}_{12}K_2=K_2R^{\tau_2}_{21}K_1R_{12},
\end{equation}
where $\tau$ acts on the corresponding tensor factor of $R$. In the vector representation $\tau$ is conjugation
by the matrix $\sigma$ interchanging $v_{N/2}$ and $v_{N/2+1}$, and $\sigma\ot\sigma$ commutes with $R$. Therefore $K$ solves~\eqref{eq:RE}
if and only if $K\sigma$ solves~\eqref{eq:tRE}. Below we write out the solutions $K$ of equation~\eqref{eq:RE}.

Note that the reflection equation can be regarded as a set of relations between the matrix entries of $K$, which generate the reflection equation algebra. This algebra can be embedded into the quantum group. In particular, from a solution of the reflection equation one constructs a coideal subalgebra with respect to which this solution is invariant (see~\cite{K}, \cite{AMS}).

Equation~\eqref{eq:RE} is quadratic and homogeneous in $K$: together with $K$, the matrix $tK$ is a solution, but the sum of two solutions is,
in general, not a solution. The reflection equation always has the trivial solution, a scalar matrix, that is, a multiple of the identity matrix.

\section{Coideal subalgebras}

\subsection{General derivation}\label{ss:gen}
We analyze the construction of a coideal subalgebra, consider the simplest examples of coideal subalgebras and show by an explicit computation an example of degeneration of the subalgebra into the whole algebra. In this section $\alpha$ is a white root and the root $\ta$ is simple (this is so in all eight cases considered below). The simplest case of a composite root is analyzed in Section~\ref{ss:comp}, and the general one in~\cite[Section~5]{AMS}.

At $q=1$ the coproduct of any element $x\in\g$ equals $x\ot1+1\ot x$, so the subalgebra generated by the elements~\eqref{eq:cl} is automatically closed under $\Delta$. For $q\neq1$ this is not so. If formula~\eqref{eq:cl} is transferred without changes, then for the element $e_\alpha+cf_{\ta}$ we obtain from~\eqref{eq:cop}
\[
\Delta(e_\alpha+cf_{\ta})=q^{h_\alpha}\ot e_\alpha+e_\alpha\ot1+c\,f_{\ta}\ot q^{-h_{\ta}}+c\,1\ot f_{\ta}.
\]
The factors $e_\alpha$ and $f_{\ta}$ are in different terms, while the subalgebra contains only their combination. For the coideal condition to hold, the second factor in every term must lie in the subalgebra. Two corrections are made.

\textbf{Step 1: the replacement $f_{\ta}\to F_{\ta}=q^{h_{\ta}}f_{\ta}$.} Since $\Delta$ is a homomorphism,
\begin{align*}
\Delta(F_{\ta})&=\Delta(q^{h_{\ta}})\Delta(f_{\ta})=(q^{h_{\ta}}\ot q^{h_{\ta}})(f_{\ta}\ot q^{-h_{\ta}}+1\ot f_{\ta})\\
&=F_{\ta}\ot1+q^{h_{\ta}}\ot F_{\ta}.
\end{align*}
When the brackets are expanded by rule~\eqref{eq:koszul}, the sign is determined by the parity of the second factor of the left bracket, that is, of $q^{h_{\ta}}$. It is even, so no sign appears, even if $f_{\ta}$ is odd.

\textbf{Step 2: the replacement $e_\alpha\to T_\alpha e_\alpha$, where $T_\alpha=q^{h_{\ta}-h_\alpha}=q^{h_{\ta}}q^{-h_\alpha}$.}
\begin{align*}
\Delta(T_\alpha e_\alpha)&=(T_\alpha\ot T_\alpha)(q^{h_\alpha}\ot e_\alpha+e_\alpha\ot1)\\
&=q^{h_{\ta}}\ot T_\alpha e_\alpha+T_\alpha e_\alpha\ot T_\alpha ,
\end{align*}
since $T_\alpha q^{h_\alpha}=q^{h_{\ta}}$. Again no sign appears, since $T_\alpha$ is even.

\textbf{Result.} Now the first factors at $T_\alpha e_\alpha$ and at $F_{\ta}$ coincide, and for
\begin{equation}\label{eq:X}
X_\alpha=T_\alpha e_\alpha+c_\alpha F_{\ta}
\end{equation}
we obtain
\begin{equation}\label{eq:DX}
\Delta(X_\alpha)=q^{h_{\ta}}\ot X_\alpha+T_\alpha e_\alpha\ot T_\alpha+c_\alpha F_{\ta}\ot1 .
\end{equation}
The second tensor factors are $X_\alpha$, $T_\alpha$ and $1$, and moreover $\Delta(T_\alpha)=T_\alpha\ot T_\alpha$. Therefore the subalgebra generated by the elements $X_\alpha$ and $T_\alpha^{\pm1}$ for all white roots $\alpha$ is a left coideal. It
is denoted by $\Uq(\kk)$ (cf.\ Theorem~5.3 in~\cite{AMS}). If $\ta=\alpha$, then $T_\alpha=1$ and
$\Delta(X_\alpha)=q^{h_\alpha}\ot X_\alpha+X_\alpha\ot1$. As $q\to1$ we have $T_\alpha\to1$, $F_{\ta}\to f_{\ta}$, and
$X_\alpha$ turns into $x_\alpha$.

\subsection{Example: $\osp(1|2)$}
Parities $(1,0,1)$, weights $\delta_1$, $0$, $-\delta_1$. The matrix $E_{12}$ has parity $|1|+|2|=1$, so $e_1$ and $f_1$ are odd, and the relation contains an anticommutator:
\begin{gather*}
e_1f_1=-E_{11}+E_{22},\qquad f_1e_1=-E_{22}+E_{33},\\
e_1f_1+f_1e_1=\diag(-1,0,1)=\frac{q^{h_1}-q^{-h_1}}{q-q^{-1}},
\end{gather*}
that is, $h_1=\diag(-1,0,1)$. The node is white, $\ta_1=\alpha_1$, $T_1=1$:
\begin{gather*}
F_1=q^{h_1}f_1=q(-E_{21}+E_{32}),\\
X_1=e_1+c_1F_1=q^{-1}E_{12}+E_{23}+c_1q(-E_{21}+E_{32}),\\
\Delta(X_1)=q^{h_1}\ot X_1+X_1\ot1 .
\end{gather*}

\subsection{Example: $\so(6)$, twisted tail}\label{ss:so6}
Here $\tau(\alpha_2)=\alpha_3$, therefore $\ta_1=\alpha_1$, $\ta_2=\alpha_3$, $\ta_3=\alpha_2$ and
\[
T_1=1,\qquad T_2=q^{h_3-h_2}=\diag(1,1,q^2,q^{-2},1,1),\qquad T_3=T_2^{-1}.
\]
Computations: $T_2e_2=E_{23}-q^{-2}E_{45}$, $T_3e_3=E_{24}-q^{-2}E_{35}$,
$F_1=q^{-1}(E_{21}-E_{65})$, $F_2=q^{-1}(E_{32}-E_{54})$, $F_3=q^{-1}(E_{42}-E_{53})$. Generators:
\begin{gather*}
X_1=e_1+c_1F_1=E_{12}-E_{56}+c_1q^{-1}(E_{21}-E_{65}),\\
X_2=T_2e_2+c_2F_3=E_{23}-q^{-2}E_{45}+c_2q^{-1}(E_{42}-E_{53}),\\
X_3=T_3e_3+c_3F_2=E_{24}-q^{-2}E_{35}+c_3q^{-1}(E_{32}-E_{54}),
\end{gather*}
and $T_2^{\pm1}$. Coproducts by formula~\eqref{eq:DX}:
\begin{gather*}
\Delta(X_1)=q^{h_1}\ot X_1+X_1\ot1,\\
\Delta(X_2)=q^{h_3}\ot X_2+T_2e_2\ot T_2+c_2F_3\ot1,\\
\Delta(X_3)=q^{h_2}\ot X_3+T_3e_3\ot T_3+c_3F_2\ot1 .
\end{gather*}
From $T_2E_{ij}T_2^{-1}=(t_i/t_j)E_{ij}$, where $t_i$ are the diagonal entries of $T_2$, we obtain $T_2X_2T_2^{-1}=q^{-2}X_2$,
$T_2X_3T_2^{-1}=q^{2}X_3$, $T_2X_1=X_1T_2$.
The same in terms of weights:
\[
T_2e_2T_2^{-1}=q^{(\alpha_3-\alpha_2,\alpha_2)}e_2=q^{-2}e_2,\qquad
T_2f_3T_2^{-1}=q^{-(\alpha_3-\alpha_2,\alpha_3)}f_3=q^{-2}f_3 .
\]
The factors coincide, so $T_2X_2T_2^{-1}$ is proportional to $X_2$. For the pair $e_2$, $f_2$ the factors would be different:
$q^{-2}$ and $q^{2}$.

Classical limit: $x_2=e_2+c_2f_3$, $x_3=e_3+c_3f_2$. Since $(\alpha_2,\alpha_3)=0$, we have $[e_2,e_3]=[f_3,f_2]=0$ and
\begin{equation}\label{eq:so6br}
[x_2,x_3]=c_3h_2-c_2h_3 .
\end{equation}
For $c_2\neq c_3$ this element together with $h_3-h_2$ gives $h_2$ and $h_3$ separately. Then $[h_2,x_2]=2e_2$ and $[h_3,x_3]=2e_3$, so $e_2$, $f_3$, $e_3$, $f_2$ get into $\kk$ separately, and from $[h_2,x_1]=-e_1+c_1f_1$ also $e_1$ and $f_1$. Hence $\kk=\g$. Therefore for $\so(6)$ the only meaningful case is $c_2=c_3$.

\subsection{Example: $\osp(2|4)$, twisted tail}
Weights of the basis vectors: $\delta_1$, $\delta_2$, $\varepsilon_1$, $-\varepsilon_1$, $-\delta_2$, $-\delta_1$. Parities
$(1,1,0,0,1,1)$. Simple roots: $\alpha_1=\delta_1-\delta_2$ is even, $\alpha_2=\delta_2-\varepsilon_1$ and $\alpha_3=\delta_2+\varepsilon_1$ are odd,
\begin{gather*}
(\alpha_1,\alpha_1)=-2,\qquad(\alpha_2,\alpha_2)=(\alpha_3,\alpha_3)=0,\\
(\alpha_1,\alpha_2)=(\alpha_1,\alpha_3)=1,\qquad(\alpha_2,\alpha_3)=-2 .
\end{gather*}
Differences from $\so(6)$: the generators $e_2$, $e_3$, $f_2$, $f_3$ are odd, $e_2^2=e_3^2=f_2^2=f_3^2=0$ (the standard
relations for odd roots with $(\alpha,\alpha)=0$, for the matrices of Section~\ref{ss:rep} they are checked by direct multiplication) and $e_2f_3+f_3e_2=0$, and $(\alpha_2,\alpha_3)=-2$ instead of $0$. The permutation $\tau$ is the same: $\tau(\alpha_2)=\alpha_3$,
$T_2=q^{h_3-h_2}=\diag(1,1,q^2,q^{-2},1,1)$, $T_3=T_2^{-1}$.

Computations: $T_2e_2=q^{-1}E_{23}+q^{-2}E_{45}$, $T_3e_3=q^{-1}E_{24}+q^{-2}E_{35}$,
$F_1=q(E_{21}-E_{65})$, $F_2=-E_{32}+qE_{54}$, $F_3=-E_{42}+qE_{53}$. Generators:
\begin{gather*}
X_1=e_1+c_1F_1=E_{12}-E_{56}+c_1q(E_{21}-E_{65}),\\
X_2=T_2e_2+c_2F_3=q^{-1}E_{23}+q^{-2}E_{45}+c_2(-E_{42}+qE_{53}),\\
X_3=T_3e_3+c_3F_2=q^{-1}E_{24}+q^{-2}E_{35}+c_3(-E_{32}+qE_{54}),
\end{gather*}
and $T_2^{\pm1}$. The element $X_1$ is even, $X_2$ and $X_3$ are odd. $X_2$ is homogeneous because $e_2$ and $f_3$ are both odd.
Coproducts:
\begin{gather*}
\Delta(X_2)=q^{h_3}\ot X_2+T_2e_2\ot T_2+c_2F_3\ot1,\\
\Delta(X_3)=q^{h_2}\ot X_3+T_3e_3\ot T_3+c_3F_2\ot1 .
\end{gather*}
Relations: $T_2X_1=X_1T_2$, $T_2X_2T_2^{-1}=q^{-2}X_2$, $T_2X_3T_2^{-1}=q^{2}X_3$ (as in $\so(6)$) and, moreover,
\[
X_2^2=X_3^2=0 .
\]
On matrices: in $X_2^2$ there are two nonzero products, $E_{45}E_{53}$ with coefficient $+c_2q^{-1}$ and $E_{42}E_{23}$ with
coefficient $-c_2q^{-1}$. Both are equal to $E_{43}$, so these terms cancel. In the algebra:
\[
X_2^2=(T_2e_2)^2+c_2^2F_3^2+c_2q^2T_2q^{h_3}(e_2f_3+f_3e_2)=0 .
\]

Classical limit: $x_2=e_2+c_2f_3$, $x_3=e_3+c_3f_2$ and, in contrast to~\eqref{eq:so6br},
\[
[x_2,x_3]=x_2x_3+x_3x_2=2E_{25}-2c_2c_3E_{52}+c_3h_2+c_2h_3 ,
\]
since $[e_2,e_3]=2E_{25}\neq0$ (the root $\alpha_2+\alpha_3=2\delta_2$). This element does not lie in the Cartan subalgebra, and the argument that led to $\kk=\g$ for $\so(6)$ does not work. That $\kk\neq\g$ for any nonzero $c_2$, $c_3$ follows from the result of Section~\ref{ss:res}: the matrix $K$ found there is not scalar at $q=1$ and commutes with $x_1$, $x_2$, $x_3$ and $h_3-h_2$.

\subsection{Generators for all eight algebras}
The generators $X_i$ in the vector representation are collected in Table~\ref{t:gens}. For the twisted diagrams $T$ is also
given.
\begin{table}[H]
\centering\small
\begin{tabular}{@{}ll@{}}
\toprule
$\so(3)$ & $X_1=E_{12}-E_{23}+c_1(E_{21}-q^{-1}E_{32})$\\[3pt]
$\so(5)$ & $X_1=E_{12}-E_{45}+c_1q^{-1}(E_{21}-E_{54})$\\
 & $X_2=E_{23}-E_{34}+c_2(E_{32}-q^{-1}E_{43})$\\[3pt]
$\so(4)$ & $X_1=E_{12}-E_{34}+c_1q^{-1}(E_{21}-E_{43})$\\
 & $X_2=E_{13}-E_{24}+c_2q^{-1}(E_{31}-E_{42})$\\[3pt]
$\so(6)$ & $T_2=q^{h_3-h_2}=\diag(1,1,q^2,q^{-2},1,1)$\\
 & $X_1=E_{12}-E_{56}+c_1q^{-1}(E_{21}-E_{65})$\\
 & $X_2=E_{23}-q^{-2}E_{45}+c_2q^{-1}(E_{42}-E_{53})$\\
 & $X_3=E_{24}-q^{-2}E_{35}+c_3q^{-1}(E_{32}-E_{54})$\\[3pt]
$\osp(1|2)$ & $X_1=q^{-1}E_{12}+E_{23}+c_1q(-E_{21}+E_{32})$\\[3pt]
$\osp(1|4)$ & $X_1=E_{12}-E_{45}+c_1q(E_{21}-E_{54})$\\
 & $X_2=q^{-1}E_{23}+E_{34}+c_2q(-E_{32}+E_{43})$\\[3pt]
$\osp(2|2)$ & $T_1=q^{h_2-h_1}=\diag(1,q^2,q^{-2},1)$\\
 & $X_1=T_1e_1+c_1F_2=q^{-1}E_{12}+q^{-2}E_{34}+c_1(-E_{31}+qE_{42})$\\
 & $X_2=T_2e_2+c_2F_1=q^{-1}E_{13}+q^{-2}E_{24}+c_2(-E_{21}+qE_{43})$\\[3pt]
$\osp(2|4)$ & $T_2=q^{h_3-h_2}=\diag(1,1,q^2,q^{-2},1,1)$\\
 & $X_1=E_{12}-E_{56}+c_1q(E_{21}-E_{65})$\\
 & $X_2=q^{-1}E_{23}+q^{-2}E_{45}+c_2(-E_{42}+qE_{53})$\\
 & $X_3=q^{-1}E_{24}+q^{-2}E_{35}+c_3(-E_{32}+qE_{54})$\\
\bottomrule
\end{tabular}
\caption{Generators $X_i$ of the coideal subalgebra $\Uq(\kk)$ (in the vector representation).}\label{t:gens}
\end{table}
For $\osp(2|2)$ here $T_2=T_1^{-1}$. For the untwisted diagrams $\Delta(X_i)=q^{h_i}\ot X_i+X_i\ot1$.

\subsection{Composite root}\label{ss:comp}
If the diagram has black nodes, the root $\ta$ may be composite. Then the element $F_{\ta}$ is built from the elements $F_\beta$ of simple roots by means of $q$-commutators. Below we analyze the simplest case, when $\ta$ is a sum of two simple roots.
For the general case see~\cite[Section~5]{AMS}. Let $A$, $B$ be elements of weights $-\beta$, $-\gamma$ such that
\[
\Delta(A)=A\ot1+q^{h_\beta}\ot A,\qquad \Delta(B)=B\ot1+q^{h_\gamma}\ot B,
\]
and $[A,B]_{q'}=AB-(-1)^{|A||B|}q'BA$. Then
\begin{multline}\label{eq:qc}
\Delta([A,B]_{q'})=[A,B]_{q'}\ot1+q^{h_\beta+h_\gamma}\ot[A,B]_{q'}\\
+(1-q'q^{-(\beta,\gamma)})\,Aq^{h_\gamma}\ot B\\
+(-1)^{|A||B|}(1-q'q^{(\beta,\gamma)})\,q^{h_\beta}B\ot A .
\end{multline}
This formula is obtained by a direct computation. The product
$\Delta(A)\Delta(B)$ is expanded by rule~\eqref{eq:koszul} into four terms:
\[
\Delta(A)\Delta(B)=AB\ot1+Aq^{h_\gamma}\ot B+(-1)^{|A||B|}q^{h_\beta}B\ot A+q^{h_\beta+h_\gamma}\ot AB ,
\]
and $\Delta(B)\Delta(A)$ is expanded in the same way. The terms of the form $(\ldots)\ot1$ and $q^{h_\beta+h_\gamma}\ot(\ldots)$ combine into the first two terms of~\eqref{eq:qc}. In the two remaining terms the elements $q^{h_\gamma}$ and $q^{h_\beta}$ are moved through $A$ and $B$ by the rules $q^{h_\gamma}A=q^{-(\beta,\gamma)}Aq^{h_\gamma}$ and $Bq^{h_\beta}=q^{(\beta,\gamma)}q^{h_\beta}B$ (they follow from the fact that $A$ and $B$ have weights $-\beta$ and $-\gamma$). This is where the factors in parentheses come from.

The first two terms of~\eqref{eq:qc} have the same form as those of $A$ and $B$ themselves. The last two terms are superfluous. The parameter $q'$ is chosen so that one of them vanishes and the remaining one has in its second factor an element of the black root (it lies in the subalgebra). If this is $A$, one takes $q'=q^{(\beta,\gamma)}$, and if $B$, then $q'=q^{-(\beta,\gamma)}$. In~\cite{AMS} the index of the $q$-commutator is understood differently:
$[F_\alpha,F_\beta]_{q^{\pm1}}=F_\alpha F_\beta-(-1)^{|F_\alpha||F_\beta|}q^{\mp(\alpha,\beta)}F_\beta F_\alpha$.

\textbf{Example: $\osp(1|4)$, the node $\alpha_1$ is black, $\alpha_2$ is white.} The reflection with respect to $\alpha_1$ interchanges
$\delta_1$ and $\delta_2$, therefore $\ta_2=\delta_1=\alpha_1+\alpha_2$ and $T_2=q^{h_{\ta_2}-h_2}=q^{h_1}$. Here
$(\alpha_1,\alpha_2)=1$, the black $F_1$ stands on the left, and one needs $q'=q$. Computations:
\begin{gather*}
F_1=q(E_{21}-E_{54}),\qquad F_2=q(-E_{32}+E_{43}),\\
F_1F_2=-q^2E_{53},\qquad F_2F_1=-q^2E_{31},\\
F_{\ta_2}=[F_1,F_2]_q=F_1F_2-qF_2F_1=q^3E_{31}-q^2E_{53}
\end{gather*}
(in the notation of~\cite{AMS} this is $[F_1,F_2]_{q^{-1}}$). From~\eqref{eq:qc}:
\[
\Delta(F_{\ta_2})=F_{\ta_2}\ot1+q^{h_1+h_2}\ot F_{\ta_2}+(1-q^2)\,q^{h_1}F_2\ot F_1 ,
\]
the superfluous term has the black $F_1$ in its second factor. Further, $q^{h_1}e_2=E_{23}+E_{34}$ and
\[
X_2=q^{h_1}e_2+c_2F_{\ta_2}=E_{23}+E_{34}+c_2(q^3E_{31}-q^2E_{53}).
\]
As in Section~\ref{ss:gen}, $\Delta(q^{h_1}e_2)=q^{h_1+h_2}\ot q^{h_1}e_2+q^{h_1}e_2\ot q^{h_1}$, therefore
\[
\Delta(X_2)=q^{h_1+h_2}\ot X_2+q^{h_1}e_2\ot q^{h_1}+c_2F_{\ta_2}\ot1+c_2(1-q^2)\,q^{h_1}F_2\ot F_1 .
\]
The second factors $X_2$, $q^{h_1}$, $1$ and $F_1=q^{h_1}f_1$ lie in the subalgebra generated by the elements $e_1$, $f_1$, $q^{\pm h_1}$,
$X_2$, therefore it is a left coideal.

\begin{propn}\label{p:coideal}
The subalgebras generated by the elements of Table~\ref{t:gens} are left coideal subalgebras. All of them, except possibly the case of $\so(6)$ with $c_2\neq c_3$, are proper.
\end{propn}
That these are left coideal subalgebras is checked by a direct computation and is essentially explained above (see~\eqref{eq:DX}). That they are proper will be proved in Section~\ref{ss:res} by exhibiting nontrivial invariants, the $K$-matrices: since the representation $V$ is irreducible, an even matrix commuting with the whole algebra $\Uq(\g)$ is scalar. For $\so(6)$ with $c_2\neq c_3$ the centralizer consists of scalar matrices, and this argument does not apply. At $q=1$ in this case $\kk=\g$ (Section~\ref{ss:so6}).

\section{$K$-matrices}\label{s:K}

\subsection{$R$-matrices}
For $\osp(1|2)$, for example, formula~\eqref{eq:R} gives the following expression for the quantum $R$-matrix:
\begin{multline*}
R=q^{-1}(E_{11}\ot E_{11}+E_{33}\ot E_{33})+q\,(E_{11}\ot E_{33}+E_{33}\ot E_{11})\\
+\sum\nolimits'E_{ii}\ot E_{jj}
+\omega\bigl(-E_{21}\ot E_{12}-q\,E_{21}\ot E_{23}\\
-(1+q)E_{31}\ot E_{13}+E_{32}\ot E_{23}+E_{32}\ot E_{12}\bigr),
\end{multline*}
where $\sum'$ is the sum over the remaining five pairs $(i,j)$.


\subsection{Finding the solutions of the reflection equation}
We solve the twisted reflection equation~\eqref{eq:tRE}: for diagrams without arcs $\tau=\mathrm{id}$, and this is
equation~\eqref{eq:RE}, while for a twisted tail its solutions are the matrices $K\sigma$, where $K$ solves~\eqref{eq:RE}
(Section~\ref{ss:RE}). Therefore in all cases below we solve equation~\eqref{eq:RE}.
The solutions of the reflection equation were found with a program in the Maple system. For each algebra the following computations were carried out.
\begin{enumerate}
\item First the generators $X_i$ (and $T$) were constructed by formula~\eqref{eq:X}, where the parameters $c_i$ are free symbols.
\item Then the linear system $KX_i=X_iK$, $KT=TK$ was solved for the entries of an even matrix $K$. Its solutions form the
\emph{centralizer}, the set of even matrices commuting with $\Uq(\kk)$. The answer was found in general form with free parameters $b_k$ (these are the free unknowns of the linear system).
\item The expression found for the centralizer was substituted into the reflection equation~\eqref{eq:RE}, after which the resulting system of quadratic equations
in $b_k$ was solved and analyzed.
\end{enumerate}
Each solution of step~3 is a line $\{tK\}$, $t\neq0$, in the centralizer. For each solution the determinant is computed. The parameters $c_i$ define a family of coideal subalgebras and, accordingly, a family of solutions. The computations were carried out with symbolic $q$ and $c_i$. It is assumed that $c_i\neq0$ and that $q$ is not a root of unity. That the matrix found commutes with the generators and satisfies the reflection equation, as well as the value of its determinant, are identities in the parameters. They hold for all nonzero $c_i$ (and, under a condition of the form $c_2=c_3$, for all values satisfying it). The dimension of the centralizer and the list of solutions are obtained by solving linear and quadratic systems with symbolic coefficients and refer to generic values of the parameters, that is, to all values except the zeros of finitely many nonzero polynomials in the parameters (for the precise definition see Section~\ref{ss:special}).

\subsection{The result without a Cartan term, formulas for $K$-matrices}\label{ss:res}
Table~\ref{t:res} gives the dimensions of the centralizers and all solutions of the reflection equation lying in the centralizer ($t$ is an arbitrary nonzero number). After the table, for each algebra a basis of the centralizer in the representation is written out. Where the dimension equals two, this basis is $1$ and $K$. For $\so(5)$, $\so(6)$ with $c_2=c_3$, $\osp(1|4)$ and $\osp(2|4)$ with $c_2=c_3$ the dimension equals three, and the basis consists of $1$, $K$ and one more matrix $K'$. The matrix $K'$ does not satisfy the reflection equation and is not a $K$-matrix. The cases $\so(3)$, $\so(4)$ and $\so(6)$ with $c_2\neq c_3$ are arranged differently and are described separately. All matrices $K$ are invertible for all nonzero $c_i$: their determinants are monomials in $q$ and $c_i$, and they are written out after the corresponding formulas. We stress that the table gives the dimensions of the centralizers in $\mathrm{End}(V)$, that is, in the representation. The coideal subalgebras $\Uq(\kk)\subset\Uq(\g)$ themselves are infinite-dimensional, and it is natural to describe them by the number of generators or by the Gelfand--Kirillov dimension.

\begin{table}[H]
\centering\small
\begin{tabular}{@{}lcl@{}}
\toprule
algebra & $\dim$ of the centralizer & solutions of the reflection equation\\
\midrule
$\so(3)$ & 3 & $t\cdot1$, $tK$\\
$\so(5)$ & 3 & $t\cdot1$, $tK$\\
$\so(4)$ & 4 & $t\cdot1$, $tX_1$, $tX_2$, $tX_1X_2$\\
$\so(6)$, $c_2=c_3$ & 3 & $t\cdot1$, $tK$\\
$\so(6)$, $c_2\neq c_3$ & 1 & $t\cdot1$\\
$\osp(1|2)$ & 2 & $t\cdot1$, $tK$\\
$\osp(1|4)$ & 3 & $t\cdot1$, $tK$\\
$\osp(2|2)$, $c_1=c_2$ & 2 & $t\cdot1$, $tK$\\
$\osp(2|2)$, $c_1\neq c_2$ & 2 & $t\cdot1$, $tK$\\
$\osp(2|4)$, $c_2=c_3$ & 3 & $t\cdot1$, $tK$\\
$\osp(2|4)$, $c_2\neq c_3$ & 2 & $t\cdot1$, $tK$\\
\bottomrule
\end{tabular}
\caption{Dimensions of the centralizers in $\mathrm{End}(V)$ and solutions of the reflection equation.}\label{t:res}
\end{table}

The centralizers and the solutions themselves are given by the following formulas.

\textbf{$\so(3)$.} Centralizer: $\spn(1,X_1,K)$ (polynomials in $X_1$),
\begin{equation}  \label{Kso(3)}
K=\begin{pmatrix}
-(1-q^{-1})c_1 & 0 & 1 \\
0 & -c_1 & 0 \\
q^{-1}c_1^{2} & 0 & 0
\end{pmatrix}.
\end{equation}
Determinant: $\det K=q^{-1}c_1^3$.

\textbf{$\so(5)$.}
\begin{gather}
K=\begin{pmatrix} \label{Kso(5)}
q^{-1}(1-q^{-1})c_1c_2 & 0 & 0 & 0 & 1 \\
0 & q^{-1}(1-q^{-1})c_1c_2 & 0 & -q^{-1}c_1 & 0 \\
0 & 0 & q^{-1}c_1c_2 & 0 & 0 \\
0 & -q^{-2}c_1c_2^{2} & 0 & 0 & 0 \\
q^{-3}c_1^{2}c_2^{2} & 0 & 0 & 0 & 0
\end{pmatrix},\\[4pt]
K'=\begin{pmatrix} \nonumber
0 & -q^{-1}c_2 & 0 & -1 & 0 \\
-q^{-2}c_1c_2 & 0 & 0 & 0 & 1 \\
0 & 0 & 0 & 0 & 0 \\
-q^{-2}c_1c_2^{2} & 0 & 0 & 0 & c_2 \\
0 & q^{-2}c_1c_2^{2} & 0 & q^{-1}c_1c_2 & 0
\end{pmatrix}.
\end{gather}
Determinant: $\det K=q^{-7}c_1^5c_2^5$.

\textbf{$\so(4)$.} Here $X_1X_2=X_2X_1$, $X_1^2=q^{-1}c_1\cdot1$, $X_2^2=q^{-1}c_2\cdot1$. The centralizer is generated by the elements $X_1$, $X_2$, and its basis is
$1$, $X_1$, $X_2$, $X_1X_2$, where
\begin{equation}  \label{Kso(4)}
X_1X_2=\begin{pmatrix}
0 & 0 & 0 & -1 \\
0 & 0 & q^{-1}c_1 & 0 \\
0 & q^{-1}c_2 & 0 & 0 \\
-q^{-2}c_1c_2 & 0 & 0 & 0
\end{pmatrix}.
\end{equation}
The reflection equation is satisfied by the generators $X_1$, $X_2$ themselves and by their product. All three matrices are invertible: $\det X_i=q^{-2}c_i^2$, $\det(X_1X_2)=q^{-4}c_1^2c_2^2$. Linear
combinations of the matrices $1$, $X_1$, $X_2$, $X_1X_2$ with more than one nonzero coefficient are not solutions.

\textbf{$\so(6)$, $c_2=c_3$.}
\begin{gather} 
K=\begin{pmatrix}  \label{Kso(6)}
(1-q^{-2})c_1c_2 & 0 & 0 & 0 & 0 & 1 \\
0 & (1-q^{-2})c_1c_2 & 0 & 0 & -q^{-1}c_1 & 0 \\
0 & 0 & c_1c_2 & 0 & 0 & 0 \\
0 & 0 & 0 & c_1c_2 & 0 & 0 \\
0 & -q^{-1}c_1c_2^{2} & 0 & 0 & 0 & 0 \\
q^{-2}c_1^{2}c_2^{2} & 0 & 0 & 0 & 0 & 0
\end{pmatrix},\\[4pt]
K'=\begin{pmatrix} \nonumber  
0 & -q^{-1}c_2 & 0 & 0 & -1 & 0 \\
-q^{-2}c_1c_2 & 0 & 0 & 0 & 0 & 1 \\
0 & 0 & 0 & 0 & 0 & 0 \\
0 & 0 & 0 & 0 & 0 & 0 \\
-q^{-1}c_1c_2^{2} & 0 & 0 & 0 & 0 & qc_2 \\
0 & q^{-1}c_1c_2^{2} & 0 & 0 & c_1c_2 & 0
\end{pmatrix}.
\end{gather}
Determinant: $\det K=q^{-4}c_1^6c_2^6$.
For $c_2\neq c_3$ the centralizer consists of scalar matrices.

\textbf{$\osp(1|2)$.}
\begin{equation} \label{Kosp(1|2)}
K=\begin{pmatrix}
-q(1+q)c_1 & 0 & 1 \\
0 & -qc_1 & 0 \\
-q^{3}c_1^{2} & 0 & 0
\end{pmatrix} .
\end{equation}
Determinant: $\det K=-q^{4}c_1^3$.

\textbf{$\osp(1|4)$.}
\begin{gather}
K=\begin{pmatrix} \label{Kosp(1|4)}
q^{2}(1+q)c_1c_2 & 0 & 0 & 0 & 1 \\
0 & q^{2}(1+q)c_1c_2 & 0 & -qc_1 & 0 \\
0 & 0 & q^{2}c_1c_2 & 0 & 0 \\
0 & q^{4}c_1c_2^{2} & 0 & 0 & 0 \\
-q^{5}c_1^{2}c_2^{2} & 0 & 0 & 0 & 0
\end{pmatrix},\\[4pt]
K'=\begin{pmatrix} \nonumber
0 & q^{2}c_2 & 0 & -1 & 0 \\
q^{3}c_1c_2 & 0 & 0 & 0 & 1 \\
0 & 0 & 0 & 0 & 0 \\
q^{4}c_1c_2^{2} & 0 & 0 & 0 & qc_2 \\
0 & -q^{4}c_1c_2^{2} & 0 & q^{2}c_1c_2 & 0
\end{pmatrix}.
\end{gather}
Determinant: $\det K=q^{12}c_1^5c_2^5$.

\textbf{$\osp(2|2)$.} The formula below gives a solution for all nonzero $c_1$, $c_2$, including $c_1=c_2$:
\begin{equation}  \label{Kosp(2|2)}
K=\begin{pmatrix}
-q^{2}(c_1+c_2) & 0 & 0 & 1 \\
0 & -q^{2}c_2 & 0 & 0 \\
0 & 0 & -q^{2}c_1 & 0 \\
-q^{4}c_1c_2 & 0 & 0 & 0
\end{pmatrix} .
\end{equation}
Determinant: $\det K=q^{8}c_1^2c_2^2$.
For $c_1\neq c_2$ the matrix $K$ is diagonalizable: its eigenvalues $-q^2c_1$ and $-q^2c_2$ have multiplicity $2$. For
$c_1=c_2$ there is a single eigenvalue, and $K$ is not diagonalizable: its Jordan blocks have sizes $2$, $1$, $1$.

\textbf{$\osp(2|4)$.} The formula below gives a solution for all nonzero $c_1$, $c_2$, $c_3$, including $c_2=c_3$:
\begin{equation}
K=\begin{pmatrix} \label{Kosp(2|4)}
q^{3}c_1(c_2+c_3) & 0 & 0 & 0 & 0 & 1 \\
0 & q^{3}c_1(c_2+c_3) & 0 & 0 & -qc_1 & 0 \\
0 & 0 & q^{3}c_1c_3 & 0 & 0 & 0 \\
0 & 0 & 0 & q^{3}c_1c_2 & 0 & 0 \\
0 & q^{5}c_1c_2c_3 & 0 & 0 & 0 & 0 \\
-q^{6}c_1^{2}c_2c_3 & 0 & 0 & 0 & 0 & 0
\end{pmatrix} .
\end{equation}
Determinant: $\det K=q^{18}c_1^6c_2^3c_3^3$.
For $c_2\neq c_3$ the matrix $K$ is diagonalizable: its eigenvalues $q^3c_1c_2$ and $q^3c_1c_3$ have multiplicity $3$.
For $c_2=c_3$ there is a single eigenvalue, and $K$ is not diagonalizable: its Jordan blocks have sizes $2$, $2$, $1$,
$1$. For $c_2=c_3$ the centralizer also contains
\begin{equation*}  
K'=\begin{pmatrix}
0 & q^{2}c_2 & 0 & 0 & -1 & 0 \\
q^{3}c_1c_2 & 0 & 0 & 0 & 0 & 1 \\
0 & 0 & 0 & 0 & 0 & 0 \\
0 & 0 & 0 & 0 & 0 & 0 \\
q^{5}c_1c_2^{2} & 0 & 0 & 0 & 0 & q^{2}c_2 \\
0 & -q^{5}c_1c_2^{2} & 0 & 0 & q^{3}c_1c_2 & 0
\end{pmatrix},\qquad K'^2=0 .
\end{equation*}

The results obtained can be formulated as the following theorem.

\begin{theorem}\label{th:main}
The matrices given by formulas~\eqref{Kso(3)}--\eqref{Kosp(2|4)}, and for $\so(4)$ the matrices $X_1$, $X_2$, $X_1X_2$, for all nonzero values of the parameters (for $\so(6)$, for $c_2=c_3$) commute with the generators of the corresponding coideal subalgebras, satisfy the reflection equation and are invertible. At generic values of the parameters there are no solutions of the reflection equation in the centralizer other than those listed in Table~\ref{t:res}. For the rows of the table with the condition $c_1=c_2$ or $c_2=c_3$ this refers to a generic point of the corresponding hypersurface. The images in the representation $V$ of the generators of the corresponding coideal subalgebras defined by the Satake diagrams in Figure~\ref{f:satake} are given in Table~\ref{t:gens}. All these subalgebras, except possibly the case of $\so(6)$ with $c_2\neq c_3$, are proper (Proposition~\ref{p:coideal}).
\end{theorem}

\section{Cartan term}\label{s:cartan}
Above, the Cartan summand 
in the element $x_\alpha$ was not considered. In this section we add it, repeat
the computations of Section~\ref{s:K} and compute the classical limit. In~\cite[Section~4]{AMS} such a summand is allowed only for an even root $\alpha$ orthogonal to the black
roots, with $\ta=\alpha=\tau(\alpha)$. Here the computations are carried out for all nodes. The main result of the paper is
Theorem~\ref{th:cartan}.

\subsection{Generators with a Cartan term}\label{ss:Xs}
For a white root $\alpha$ we set
\begin{equation}\label{eq:Xs}
X^{(s)}_\alpha=T_\alpha e_\alpha+c_\alpha F_{\ta}+s\,(q^{h_{\ta}}-1),\qquad s\in\mathbb C .
\end{equation}
Since $\Delta(q^{h}-1)=q^{h}\ot(q^{h}-1)+(q^{h}-1)\ot1$, from~\eqref{eq:DX} we obtain
\[
\Delta(X^{(s)}_\alpha)=q^{h_{\ta}}\ot X^{(s)}_\alpha+T_\alpha e_\alpha\ot T_\alpha
+\bigl(c_\alpha F_{\ta}+s\,(q^{h_{\ta}}-1)\bigr)\ot1 ,
\]
therefore the subalgebra generated by the elements~\eqref{eq:Xs} and $T_\alpha^{\pm1}$ is a left coideal.

As $q\to1$, $(q-1)s\to\sigma$, the element $X^{(s)}_\alpha$ turns into $x_\alpha+\sigma h_{\ta}$. The subalgebra of $\g$
generated by the elements $x_\beta$ and $h_{\tilde\beta}-h_\beta$, where at one or several nodes $x_\alpha$
is replaced by $x_\alpha+\sigma_\alpha h_{\ta}$, will be denoted by $\kk_\sigma$. At $\sigma=0$ it turns into $\kk$.

\subsection{$K$-matrices}\label{ss:Ks}
The computations follow the scheme of Section~\ref{s:K}: one of the generators $X_\alpha$ is replaced by $X^{(s)}_\alpha$ (for
$\so(4)$, both of them, with independent $s_1$ and $s_2$. A Cartan term on one node only is obtained for $s_2=0$ or $s_1=0$, and the formulas of this section for $\so(4)$ remain valid in these cases), and $s$ is one more free symbol. The results are collected in
Table~\ref{t:resS} ($t$ is an arbitrary nonzero number). The nodes $\alpha_3$ in $\so(6)$, $\osp(2|4)$ and $\alpha_2$ in $\osp(2|2)$ are mapped to the ones listed by
$\tau$, and the computation for them gives the same results.

\begin{table}[H]
\centering\small
\setlength{\tabcolsep}{4pt}
\begin{tabular}{@{}llcl@{}}
\toprule
(super)algebra & Cartan term on & \hd{$\dim$ of the\\ centralizer} & \hd[l]{solutions of the\\ reflection equation}\\
\midrule
$\so(3)$ & $\alpha_1$ & 3 & $t\cdot1$, $tK(s)$\\
$\so(5)$ & $\alpha_1$ & 3 & $t\cdot1$, $tK(s)$\\
$\so(5)$ & $\alpha_2$ & 1 & $t\cdot1$\\
$\so(4)$ & $\alpha_1$, $\alpha_2$, $\alpha_1$ and $\alpha_2$ & 4 & $t\cdot1$, $tK_1$, $tK_2$, $tK_1K_2$\\
$\so(6)$, $c_2=c_3$ & $\alpha_1$ & 3 & $t\cdot1$, $tK(s)$\\
$\so(6)$, $c_2=c_3$ & $\alpha_2$ & 1 & $t\cdot1$\\
$\osp(1|2)$ & $\alpha_1$ & 1 & $t\cdot1$\\
$\osp(1|4)$ & $\alpha_1$ & 3 & $t\cdot1$, $tK(s)$\\
$\osp(1|4)$ & $\alpha_2$ & 1 & $t\cdot1$\\
$\osp(2|2)$ & $\alpha_1$ & 1 & $t\cdot1$\\
$\osp(2|4)$, $c_2=c_3$ & $\alpha_1$ & 3 & $t\cdot1$, $tK(s)$\\
$\osp(2|4)$, $c_2\neq c_3$ & $\alpha_1$ & 1 & $t\cdot1$\\
$\osp(2|4)$, $c_2=c_3$ & $\alpha_2$ & 1 & $t\cdot1$\\
\bottomrule
\end{tabular}
\caption{Cartan term: dimensions of the centralizers and solutions of the reflection equation.}\label{t:resS}
\end{table}

Below $K(s)$ denotes the nonscalar solution with a Cartan term. At $s=0$ it coincides with the matrix $K$
from Section~\ref{ss:res}.

\textbf{$\so(3)$}, term on $\alpha_1$:
\begin{align*}
K(s)&=\frac{c_1}{q}\cdot1-(1-q^{-1})\,s\,X^{(s)}_1-\bigl(X^{(s)}_1\bigr)^2\\
&=\begin{pmatrix}
-(1-q^{-1})c_1-(q-1)\omega s^2 & -\omega s & 1 \\
-c_1\omega s & -c_1 & 0 \\
q^{-1}c_1^{2} & 0 & 0
\end{pmatrix}.
\end{align*}

\textbf{$\so(4)$}, terms $s_1$, $s_2$ on both nodes. Here $X^{(s_1)}_1X^{(s_2)}_2=X^{(s_2)}_2X^{(s_1)}_1$, and
the solutions are $K_1$, $K_2$ and $K_1K_2$, where $K_i=X^{(s_i)}_i+(1-q^{-1})\,s_i\cdot1$:
\[
K_1=\begin{pmatrix}
\omega s_1 & 1 & 0 & 0 \\
q^{-1}c_1 & 0 & 0 & 0 \\
0 & 0 & \omega s_1 & -1 \\
0 & 0 & -q^{-1}c_1 & 0
\end{pmatrix},\qquad
K_2=\begin{pmatrix}
\omega s_2 & 0 & 1 & 0 \\
0 & \omega s_2 & 0 & -1 \\
q^{-1}c_2 & 0 & 0 & 0 \\
0 & -q^{-1}c_2 & 0 & 0
\end{pmatrix}.
\]

\textbf{$\so(5)$}, term on $\alpha_1$:
{\setlength{\arraycolsep}{2.5pt}
\[
K(s)=\begin{pmatrix}
q^{-1}c_2\bigl((1-q^{-1})c_1-\omega^2s^2\bigr) & -q^{-1}c_2\omega s & 0 & -\omega s & 1 \\
-q^{-2}c_1c_2\omega s & q^{-1}(1-q^{-1})c_1c_2 & 0 & -q^{-1}c_1 & 0 \\
0 & 0 & q^{-1}c_1c_2 & 0 & 0 \\
-q^{-2}c_1c_2^{2}\omega s & -q^{-2}c_1c_2^{2} & 0 & 0 & 0 \\
q^{-3}c_1^{2}c_2^{2} & 0 & 0 & 0 & 0
\end{pmatrix}.
\]
}

\textbf{$\so(6)$, $c_2=c_3$}, term on $\alpha_1$:
{\setlength{\arraycolsep}{4pt}
\[
K(s)=\begin{pmatrix}
c_2\bigl((1-q^{-2})c_1-q^{-1}\omega^2s^2\bigr) & -q^{-1}c_2\omega s & 0 & 0 & -\omega s & 1 \\
-q^{-2}c_1c_2\omega s & (1-q^{-2})c_1c_2 & 0 & 0 & -q^{-1}c_1 & 0 \\
0 & 0 & c_1c_2 & 0 & 0 & 0 \\
0 & 0 & 0 & c_1c_2 & 0 & 0 \\
-q^{-1}c_1c_2^{2}\omega s & -q^{-1}c_1c_2^{2} & 0 & 0 & 0 & 0 \\
q^{-2}c_1^{2}c_2^{2} & 0 & 0 & 0 & 0 & 0
\end{pmatrix}.
\]
}

\textbf{$\osp(1|4)$}, term on $\alpha_1$:
\[
K(s)=\begin{pmatrix}
q^{2}c_2\bigl((1+q)c_1+\omega^2s^2\bigr) & -q^{2}c_2\omega s & 0 & \omega s & 1 \\
-q^{3}c_1c_2\omega s & q^{2}(1+q)c_1c_2 & 0 & -qc_1 & 0 \\
0 & 0 & q^{2}c_1c_2 & 0 & 0 \\
-q^{4}c_1c_2^{2}\omega s & q^{4}c_1c_2^{2} & 0 & 0 & 0 \\
-q^{5}c_1^{2}c_2^{2} & 0 & 0 & 0 & 0
\end{pmatrix}.
\]

\textbf{$\osp(2|4)$, $c_2=c_3$}, term on $\alpha_1$:
\[
K(s)=\begin{pmatrix}
q^{2}c_2\bigl(2qc_1+\omega^2s^2\bigr) & -q^{2}c_2\omega s & 0 & 0 & \omega s & 1 \\
-q^{3}c_1c_2\omega s & 2q^{3}c_1c_2 & 0 & 0 & -qc_1 & 0 \\
0 & 0 & q^{3}c_1c_2 & 0 & 0 & 0 \\
0 & 0 & 0 & q^{3}c_1c_2 & 0 & 0 \\
-q^{5}c_1c_2^{2}\omega s & q^{5}c_1c_2^{2} & 0 & 0 & 0 & 0 \\
-q^{6}c_1^{2}c_2^{2} & 0 & 0 & 0 & 0 & 0
\end{pmatrix}.
\]

Determinants: $\det K(s)=\det K$, that is, $c_1^3q^{-1}$, $c_1^5c_2^5q^{-7}$, $c_1^6c_2^6q^{-4}$, $q^{12}c_1^5c_2^5$,
$q^{18}c_1^6c_2^6$ for $\so(3)$, $\so(5)$, $\so(6)$, $\osp(1|4)$, $\osp(2|4)$ respectively, and for $\so(4)$
$\det K_i=c_i^2q^{-2}$ and $\det(K_1K_2)=c_1^2c_2^2q^{-4}$. The determinants do not depend on $s$, $s_i$, therefore all these matrices are invertible for all values of $s$, $s_i$ and nonzero $c_i$.

\subsection{Special values of the parameters}\label{ss:special}
Denote by $A=\mathrm{End}_{\mathcal B}(V)$ the centralizer: it is a unital associative algebra consisting of
even matrices commuting with $\mathcal B$. Denote by $\operatorname{rad}A$ the radical, that is, the largest ideal of $A$ consisting of nilpotent elements (the sum of two such ideals again consists of nilpotent elements, so the largest ideal exists). In all the cases under consideration $A$ is commutative (this was checked
by a computation at a generic point of each of the parameter sets considered below). For a commutative $A$ the radical consists of all its nilpotent elements.

\textbf{Definition.} A statement holds \emph{at a generic point} (in general position) of an irreducible set $Y$ in the parameter space if it holds
everywhere on $Y$ outside some proper (Zariski) closed subset. A computation at a generic point is a
computation over the field of rational functions on $Y$. Let $d$ and $r$ be the values of $\dim A$ and $\dim\operatorname{rad}A$ at
a generic point of $Y$. An irreducible hypersurface $Y'\subset Y$ is called a set of \emph{special values of the first type}
if $\dim A>d$ at a generic point of $Y'$, and a set of \emph{special values of the second type} if
$\dim A=d$ and $\dim\operatorname{rad}A>r$ at a generic point of $Y'$. \\

First $Y$ is the whole parameter space (at its generic point $\operatorname{rad}A=0$ in all cases). Then the same definition is applied inside
the sets found. The parameters $c_i$ and the parameters
of the Cartan terms are assumed to be nonzero here: the zero value of a parameter corresponds to the case
without a Cartan term at this node, and such cases are considered separately. Everywhere it is assumed that $q$ is not a root of unity.

\textbf{Lemma.} Let $A\subset\mathrm{End}(V)$ be a subalgebra over $\mathbb C$ containing the identity, and let
$b_1,\dots,b_d$ be its basis. Then $\operatorname{rad}A=\{a\in A:\ \operatorname{tr}(ab)=0$ for all $b\in A\}$. In
particular, $\dim\operatorname{rad}A=d-\operatorname{rank}\bigl(\operatorname{tr}(b_ib_j)\bigr)$, and $\operatorname{rad}A\neq0$
if and only if $\det\bigl(\operatorname{tr}(b_ib_j)\bigr)=0$.

\textit{Proof.} Denote by $I$ the set on the right-hand side. It is a subspace, since
the condition $\operatorname{tr}(ab)=0$ is linear in $a$.

Let $a\in\operatorname{rad}A$ and $b\in A$. Since $\operatorname{rad}A$ is an ideal, $ab\in\operatorname{rad}A$, therefore
$ab$ is nilpotent. All eigenvalues of a nilpotent matrix are zero, and $\operatorname{tr}(ab)=0$. Hence
$\operatorname{rad}A\subset I$.

Conversely, $I$ is a two-sided ideal: if $a\in I$ and $c\in A$, then by the cyclicity of the trace
$\operatorname{tr}((ca)b)=\operatorname{tr}(a(bc))=0$ and $\operatorname{tr}((ac)b)=\operatorname{tr}(a(cb))=0$ for all
$b\in A$, since $bc,cb\in A$. Let us show that the elements of $I$ are nilpotent. Let $a\in I$. Since $A$ contains
the identity, $a^{k-1}\in A$ for all $k\geqslant1$, and $p_k:=\operatorname{tr}(a^k)=\operatorname{tr}(a\cdot a^{k-1})=0$.
If $\lambda_1,\dots,\lambda_N$ are the eigenvalues of $a$ (with multiplicities), then $p_k=\sum_i\lambda_i^k$, and Newton's
identities
\[
k\,e_k=\sum_{i=1}^{k}(-1)^{i-1}e_{k-i}\,p_i,\qquad e_0=1,
\]
where $e_k$ are the elementary symmetric polynomials in $\lambda_1,\dots,\lambda_N$, give by induction
$e_1=\dots=e_N=0$ (over $\mathbb C$ one can divide by $k$). Therefore the characteristic polynomial of $a$ equals
$\lambda^N$, and by the Cayley--Hamilton theorem $a^N=0$. Thus $I$ is an ideal of nilpotent elements, and
$I\subset\operatorname{rad}A$, since $\operatorname{rad}A$ is the largest such ideal.

Finally, an element $a=\sum_ix_ib_i$ lies in $I$ if and only if $\sum_ix_i\operatorname{tr}(b_ib_j)=0$
for all $j$. Therefore $\operatorname{rad}A=I$ is the kernel of the matrix $\bigl(\operatorname{tr}(b_ib_j)\bigr)$, whence
$\dim\operatorname{rad}A=d-\operatorname{rank}\bigl(\operatorname{tr}(b_ib_j)\bigr)$. $\square$

The special values are found by a Maple program. At a generic point
of the set $Y$ one finds a basis $b_1,\dots,b_d$ of the centralizer and $r=d-\operatorname{rank}\bigl(\operatorname{tr}(b_ib_j)\bigr)$
(by the lemma). The candidates for special values are the irreducible factors of the following expressions:
\begin{enumerate}
\item The minor of maximal order of the matrix of the linear system $KX_i=X_iK$, $KT=TK$ (the last pivot of Bareiss elimination). Outside its zeros the rank of the system, and hence $\dim A$, do not change.
\item The numerator and the denominator of the minor of maximal order of $\bigl(\operatorname{tr}(b_ib_j)\bigr)$. Outside their zeros $\dim\operatorname{rad}A$ does not change.
\item The denominators of the basis elements $b_k$. Outside their zeros the basis is defined.
\end{enumerate}
Factors equal to the parameters $c_i$, $s$, $s_i$ themselves are discarded, as are the factors depending only on $q$:
these are $q-1$, $q+1$, $q^2+1$, $q^2+q+1$, and their roots are roots of unity. Outside the zeros of the candidates, $\dim A$ and
$\dim\operatorname{rad}A$ are the same as at a generic point of $Y$. For each candidate $f$, a parameter entering $f$
linearly is expressed from the equation $f=0$ (such a parameter was found for all candidates). At a generic point of the resulting
hypersurface $\{f=0\}\subset Y$ one finds $\dim A$ and $\dim\operatorname{rad}A$, and the whole procedure
is repeated for it, whether these numbers have changed or not. Each step decreases the number of free parameters by one, therefore
the procedure is finite, and it finds all sets of special values, including those inside the sets found.
The procedure was carried out for all sets of nodes of all eight algebras, including the sets without Cartan terms. The parameters
$c_i$ are independent here, and the conditions $c_2=c_3$ for $\so(6)$, $\osp(2|4)$ and $c_1=c_2$ for $\osp(2|2)$ are obtained
as special values.

\textbf{Special values of the first type} are $c_2=c_3$ for $\so(6)$ and $\osp(2|4)$ (without a Cartan term on $\alpha_1$ and with it)
and the hypersurface~\eqref{eq:sp} below for $\so(5)$ with terms on both nodes. There are no others, including inside the sets
of special values. With terms on several nodes at once, the centralizer is nonscalar only for $\so(4)$
(Section~\ref{ss:Ks}) and under condition~\eqref{eq:sp}. For $\osp(2|4)$ without a Cartan term at $c_2=c_3$, the
radical changes together with the dimension: $\dim A=3$, $\dim\operatorname{rad}A=2$. Here $V$ is indecomposable over $\mathcal B$ (in the notation of Table~\ref{t:typeII},
$(2|4)_{2,2,1,1}$), and the solutions of the reflection equation are still $t\cdot1$ and $tK$ from Section~\ref{ss:res}.

\textbf{$\so(5)$, terms on both nodes.} We consider the generators
\[
X^{(s)}_1=X_1+s\,(q^{h_1}-1),\qquad X^{(s_2)}_2=X_2+s_2\,(q^{h_2}-1).
\]
For generic values of the parameters the centralizer consists of scalar matrices. At a generic point of the hypersurface
\begin{equation}\label{eq:sp}
c_1=-(q-1)^2s^2
\end{equation}
it is three-dimensional, and the solutions of the reflection equation in it are $t\cdot1$ and $tK(s,s_2)$, where
{\setlength{\arraycolsep}{4pt}
\[
K(s,s_2)=\begin{pmatrix}
\mu-q^{-1}c_2\omega^2s^2 & -q^{-1}c_2\omega s & q^{-1}(q-1)\omega ss_2 & -\omega s & 1 \\
-q^{-2}c_1c_2\omega s & \mu & q^{-1}c_1\omega s_2 & -q^{-1}c_1 & 0 \\
q^{-2}(q-1)c_1c_2\omega ss_2 & q^{-1}c_1c_2\omega s_2 & q^{-1}c_1c_2 & 0 & 0 \\
-q^{-2}c_1c_2^{2}\omega s & -q^{-2}c_1c_2^{2} & 0 & 0 & 0 \\
q^{-3}c_1^{2}c_2^{2} & 0 & 0 & 0 & 0
\end{pmatrix},
\]
}
$\mu=q^{-1}(1-q^{-1})c_1c_2+q^{-1}(q-1)c_1\omega s_2^2$, and everywhere $c_1$ is related to $s$ by condition~\eqref{eq:sp}. For $s_2=0$
the matrix $K(s,s_2)$ coincides with $K(s)$ from Section~\ref{ss:Ks}.
The determinant is the same as that of $K$: $\det K(s,s_2)=c_1^5c_2^5q^{-7}$. As $q\to1$, $(q-1)s\to\sigma$,
condition~\eqref{eq:sp} turns into $c_1=-\sigma^2$.

\textbf{Remark.} Under condition~\eqref{eq:sp} the generator $X^{(s)}_1$ can be written as
\[
X^{(s)}_1=e_1-a^2q^{h_1}f_1-a\,\frac{q^{h_1}-1}{q-1},\qquad a=-(q-1)s,
\]
that is, it has the form of the generators of the quantum electrical algebra from~\cite[Section~3.1]{AMS2} (there they are considered for
$\mathfrak{sl}(N)$ and the twisted reflection equation). Its classical limit
$x_1+\sigma h_1=e_1-\sigma^2f_1+\sigma[e_1,f_1]$ has the form of the generators $e_\beta-a_\beta^2f_\beta-a_\beta[e_\beta,f_\beta]$
of the electrical Lie algebra from~\cite[Section~3.1]{AMS2} with $a_\beta=-\sigma$ (on electrical Lie algebras
see~\cite{LP,BGG}).

Inside~\eqref{eq:sp} there are no special values of the first type. Special
values of the second type inside~\eqref{eq:sp} are given in Table~\ref{t:typeII}.

\textbf{Eigenvalues of the $K$-matrices under condition~\eqref{eq:sp}.} Any matrix $K\in A$ acts on each
indecomposable summand $p_iV$ as $\lambda_i\cdot1$ plus a nilpotent matrix (see the paragraph ``Decomposition of $V$'' below).
Therefore the eigenvalues of $K$ and their multiplicities are related to the decomposition of $V$ over $\mathcal B$. Under condition~\eqref{eq:sp} and
for generic $s_2$ we have $V=1+2+2$ and $\operatorname{rad}A=0$ (Table~\ref{t:dec}). The matrix $K(s,s_2)$ is diagonalizable and
has three eigenvalues of multiplicities $2$, $2$ and $1$, one on each summand:
\[
K(s,s_2)\sim\begin{pmatrix}
\lambda_1&&&&\\ &\lambda_1&&&\\ &&\lambda_2&&\\ &&&\lambda_2&\\ &&&&\lambda_3
\end{pmatrix},
\]
where $\lambda_3=q^{-3}c_1c_2$, and $\lambda_1$, $\lambda_2$ are the roots of $t^2-\bigl(2q^{-1}c_1c_2+q^{-1}(q-1)\omega c_1s_2^2\bigr)t+q^{-2}c_1^2c_2^2=0$.

 For generic values of the parameters
and a Cartan term only on $\alpha_1$ we have $V=3+1+1$. The matrix $K(s)$ has three eigenvalues of multiplicities $3$, $1$ and $1$, while for $K$
from Section~\ref{ss:res} the values on the two one-dimensional summands coincide, and the multiplicities are $3$ and $2$.

For $s_2=0$, that is, for $\so(5)$ with a Cartan term only on $\alpha_1$, condition~\eqref{eq:sp} defines special
values of the second type (Table~\ref{t:typeII}): $V=1+4_{2,1,1}$ and $\dim\operatorname{rad}A=1$. The matrix $K(s)$ is not
diagonalizable here. It has a simple eigenvalue on the one-dimensional summand and an eigenvalue of multiplicity $4$ with
Jordan blocks of sizes $2$, $1$, $1$ on the four-dimensional one:
\[
K(s)\sim\begin{pmatrix}
\lambda&1&&&\\ &\lambda&&&\\ &&\lambda&&\\ &&&\lambda&\\ &&&&\lambda_3
\end{pmatrix},\qquad \lambda=q^{-1}c_1c_2,\quad \lambda_3=q^{-3}c_1c_2 .
\]

\textbf{Classical limit under condition~\eqref{eq:sp}.} As $q\to1$, condition~\eqref{eq:sp} turns into $c_1=-\sigma^2$.
The subalgebra $\kk_\sigma$ with Cartan terms on both nodes at $c_1=-\sigma^2$ is isomorphic to $\mathfrak{gl}(2)$, and $V$ decomposes
over it into a sum of subspaces of dimensions $2$, $2$ and $1$. For generic values of the parameters, the subalgebra $\kk$ and the subalgebra
$\kk_\sigma$ with a term on $\alpha_1$ are isomorphic to $\so(3)\oplus\so(2)\cong\mathfrak{gl}(2)$, but $V$ decomposes over them
into a sum of subspaces of dimensions $3$, $1$ and $1$. For the subalgebra $\kk_\sigma$ with a term only on $\alpha_1$ at
$c_1=-\sigma^2$, the derived algebra is two-dimensional and abelian.

\textbf{Decomposition of $V$.} Let $b_1,\dots,b_d$ be a basis of $A$ at a generic point of the set under consideration and
$a=t_1b_1+\dots+t_db_d$, where $t_j$ are new independent symbols (a generic element of $A$). The decomposition of $V$ is found as follows.
\begin{enumerate}
\item Decompositions $V=V_1\oplus\dots\oplus V_n$ into a direct sum of submodules over $\mathcal B$ are in one-to-one
correspondence with sets of orthogonal idempotents $p_1,\dots,p_n\in A$ with sum $1$: $p_i$ is the projection onto $V_i$ along
the other summands. The summand $V_i$ is indecomposable if and only if $p_i$ cannot be represented as a sum of
two orthogonal nonzero idempotents from $A$.
\item Since $A$ is commutative, $A/\operatorname{rad}A\cong\mathbb C^n$ by the Artin--Wedderburn theorem, where $n=\dim A-\dim\operatorname{rad}A$.
The primitive idempotents of $A/\operatorname{rad}A$ lift uniquely to orthogonal idempotents
$p_1,\dots,p_n\in A$ with sum $1$, and each algebra $p_iA$ is local: $p_iA=\mathbb C\,p_i\oplus\operatorname{rad}(p_iA)$.
Therefore $V=\bigoplus_ip_iV$ is a decomposition into indecomposable submodules, and there are no other such decompositions: in a commutative
algebra a complete system of orthogonal primitive idempotents is unique. (For a noncommutative $A$ the decomposition is
unique only up to isomorphism of the summands, by the Krull--Schmidt theorem, see~\cite[\S3.8]{EGH}.)
\item A generic element has the form $a=\sum_i(\lambda_ip_i+\nu_i)$ with pairwise distinct $\lambda_i$ and nilpotent
$\nu_i\in p_iA$. On $p_iV$ the matrix $a-\lambda_i$ is nilpotent, and on $p_jV$ with $j\neq i$ it is invertible. Therefore the summands
$p_iV$ are exactly the generalized eigenspaces of $a$, and there are exactly $n$ of them.
\item Any matrix $K\in A$ also has the form $\sum_i(\lambda_i(K)p_i+\nu_i(K))$, but the numbers $\lambda_i(K)$ on different
summands may coincide, and the nilpotent parts may be smaller than those of a generic element. If
$\operatorname{rad}A=0$, then $A\cong\mathbb C^{\dim A}$, all $\nu_i$ are zero, every
matrix $K\in A$ is diagonalizable, and the summands are pairwise nonisomorphic, that is, $V$ decomposes over $\mathcal B$ without
multiplicities.
\end{enumerate}

The summands are found as follows. The characteristic polynomial of $a$ is factored into irreducible factors over the field
of rational functions in the parameters and $t_j$: $\chi_a=\prod_kg_k^{m_k}$. Each root $\lambda$ of a factor $g_k$ gives
the summand $\ker(a-\lambda)^{m_k}$ of dimension $m_k$. If $\deg g_k=\ell>1$, then the $\ell$ roots are conjugate, and their
summands are arranged in the same way. The Jordan blocks of $a$ on these summands are found from the ranks $R_j=\operatorname{rank}g_k(a)^j$:
for each root of $g_k$ the number of blocks of size $j$ equals $(R_{j-1}-2R_j+R_{j+1})/\ell$. Indeed, on the summand
corresponding to a root $\lambda$ of $g_k$, the matrix $g_k(a)$ differs from $a-\lambda$ by an invertible factor, and on the summands corresponding to the other factors $g_{k'}$ the matrix $g_k(a)$ is invertible. Each of the $\ell$ summands
corresponding to the roots of $g_k$ makes the same contribution, therefore $R_j=\dim V-\ell m_k+\ell\rho_j$, where
$\rho_j=\operatorname{rank}(a-\lambda)^j$ on the summand, and $\rho_{j-1}-\rho_j$ is the number of blocks of size at least $j$.
For superalgebras $g_k(a)^{m_k}$ is an even matrix. Its kernel is the sum of $\ell$ summands, and the even part of each of them
has dimension $(n_0-\operatorname{rank}P_0)/\ell$, where $n_0=\dim V_{\bar0}$ and $P_0$ is the block of the matrix
$g_k(a)^{m_k}$ on $V_{\bar0}$. A summand of superdimension $(m_0|m_1)$ (for $\so(n)$, of dimension $m$) on which
$a$ has Jordan blocks of sizes $j_1\geqslant j_2\geqslant\dots$, not all equal to $1$, is written as
$(m_0|m_1)_{j_1,j_2,\dots}$ (respectively $m_{j_1,j_2,\dots}$). For example, $4_{2,1,1}$ is a summand of dimension $4$
with blocks of sizes $2$, $1$, $1$. In each case the number of summands found coincided with
$\dim A-\dim\operatorname{rad}A$. The superdimensions of the summands at a generic point are given in Table~\ref{t:dec}:

\begin{table}[H]
\centering\small
\begin{tabular}{@{}llcl@{}}
\toprule
algebra & Cartan term on & $\dim A$ & decomposition of $V$\\
\midrule
$\so(3)$ & none, $\alpha_1$ & 3 & $1+1+1$\\
$\so(4)$ & none, $\alpha_1$, $\alpha_2$, $\alpha_1$ and $\alpha_2$ & 4 & $1+1+1+1$\\
$\so(5)$ & none, $\alpha_1$ & 3 & $3+1+1$\\
$\so(5)$, \eqref{eq:sp} & $\alpha_1$ and $\alpha_2$ & 3 & $1+2+2$\\
$\so(6)$, $c_2=c_3$ & none, $\alpha_1$ & 3 & $4+1+1$\\
$\osp(1|2)$ & none & 2 & $(1|1)+(0|1)$\\
$\osp(1|4)$ & none, $\alpha_1$ & 3 & $(1|2)+(0|1)+(0|1)$\\
$\osp(2|2)$, $c_1\neq c_2$ & none & 2 & $(1|1)+(1|1)$\\
$\osp(2|4)$, $c_2\neq c_3$ & none & 2 & $(1|2)+(1|2)$\\
$\osp(2|4)$, $c_2=c_3$ & $\alpha_1$ & 3 & $(2|2)+(0|1)+(0|1)$\\
\bottomrule
\end{tabular}
\caption{Decomposition of $V$ over $\mathcal B$ for generic values of the parameters.}\label{t:dec}
\end{table}

In the remaining cases the centralizer at a generic point consists of scalar matrices, and $V$ is indecomposable over $\mathcal B$.

\textbf{Special values of the second type} are listed in Table~\ref{t:typeII}. On each hypersurface of the table
$\dim\operatorname{rad}A=1$, except for $\so(4)$, where $\dim\operatorname{rad}A=2$. Here two summands of the decomposition of $V$
are replaced by one summand on which a generic element of $A$ is not diagonalizable (for $\so(4)$ this happens with two
pairs of summands). For $\so(4)$ the hypersurface $4qc_i=-(q^2-1)^2s_i^2$ is the same for a term only on $\alpha_i$
and for terms on both nodes. Inside the hypersurfaces of the table, new special values occur only for $\so(4)$ with
terms on both nodes: at the intersection of the two hypersurfaces $\dim A=4$, $\dim\operatorname{rad}A=3$, and $V$
is indecomposable (the last row for $\so(4)$ in Table~\ref{t:typeII}). The last
column gives the limit of the condition as $q\to1$, $(q-1)s\to\sigma$, $(q-1)s_i\to\sigma_i$. For
$\so(3)$ the conditions of the table are related to the characteristic polynomial of the matrix $K(s)$:
\[
\det\bigl(K(s)-t\bigr)=-\bigl(t-q^{-1}c_1\bigr)\bigl(t^2+(2c_1+(q-1)\omega s^2)\,t+c_1^2\bigr).
\]
The first condition is the vanishing of the discriminant $(q-1)\omega s^2\bigl(4c_1+(q-1)\omega s^2\bigr)$
of the quadratic factor, and the second is that this factor has the root $q^{-1}c_1$.

\begin{table}[H]
\centering\small
\setlength{\tabcolsep}{4pt}
\begin{tabular}{@{}lllc@{}}
\toprule
algebra, Cartan term on & condition & decomposition of $V$ & limit\\
\midrule
$\so(3)$, $\alpha_1$ & $4qc_1=-(q-1)^2(q+1)s^2$ & $1+2_{2}$ & $c_1=-\sigma^2/2$\\
 & $(q+1)c_1=-(q-1)^2s^2$ & $1+2_{2}$ & $c_1=-\sigma^2/2$\\
$\so(4)$, $\alpha_i$ or $\alpha_1,\alpha_2$ & $4qc_i=-(q^2-1)^2s_i^2$ & $2_{2}+2_{2}$ & $c_i=-\sigma_i^2$\\
 & the same for $i=1$ and $i=2$ & $4_{3,1}$ & $c_i=-\sigma_i^2$\\
$\so(5)$, $\alpha_1$ & $4qc_1=-(q^2-1)^2s^2$ & $3+2_{2}$ & $c_1=-\sigma^2$\\
 & $c_1=-(q-1)^2s^2$ & $1+4_{2,1,1}$ & $c_1=-\sigma^2$\\
$\so(5)$, \eqref{eq:sp}, $\alpha_1$ and $\alpha_2$ & $4qc_2=-(q-1)^2(q+1)s_2^2$ & $1+4_{2,2}$ & $c_2=-\sigma_2^2/2$\\
 & $(q+1)c_2=qs_2^2$ & $2+3_{2,1}$ & none\\
$\so(6)$, $c_2=c_3$, $\alpha_1$ & $4qc_1=-(q^2-1)^2s^2$ & $4+2_{2}$ & $c_1=-\sigma^2$\\
 & $(q^2+1)^2c_1=-q(q^2-1)^2s^2$ & $1+5_{2,1,1,1}$ & $c_1=-\sigma^2$\\
$\osp(1|4)$, $\alpha_1$ & $4q^3c_1=-(q^2-1)^2s^2$ & $(1|2)+(0|2)_{2}$ & $c_1=-\sigma^2$\\
 & $q^2c_1=(q+1)^2s^2$ & $(0|1)+(1|3)_{2,1,1}$ & none\\
$\osp(2|2)$, none & $c_1=c_2$ & $(2|2)_{2,1,1}$ & ---\\
$\osp(2|4)$, $c_2=c_3$, $\alpha_1$ & $4q^3c_1=-(q^2-1)^2s^2$ & $(2|2)+(0|2)_{2}$ & $c_1=-\sigma^2$\\
\bottomrule
\end{tabular}
\caption{Special values of the second type and the decomposition of $V$ at a generic point of the corresponding set. The second row for
$\so(4)$ is the intersection of the two hypersurfaces with terms on both nodes (codimension $2$).}\label{t:typeII}
\end{table}

In the two cases marked ``none'', the value of $c_i$ grows without bound under the indicated limit transition.

\textbf{$K$-matrices on the sets of special values.} At a generic point of each set of special values the reflection
equation was solved anew following the scheme of Section~\ref{s:K}. On each set of the second type the solutions are the same lines as at
a generic point of the ambient set, with the condition substituted: $t\cdot1$ and $tK(s)$ (for $\so(5)$ inside~\eqref{eq:sp}
these are $t\cdot1$ and $tK(s,s_2)$, for $\osp(2|2)$ at $c_1=c_2$ they are $t\cdot1$ and $tK$ from Section~\ref{ss:res}, and for $\so(4)$ they are $t\cdot1$, $tK_1$, $tK_2$,
$tK_1K_2$). There are no other solutions. Their determinants are the same monomials with the condition substituted, therefore these matrices
are invertible. The Jordan blocks of the matrix $K(s)$ (respectively $K(s,s_2)$ or $K$, and for $\so(4)$ of the matrix $K_1K_2$) at
a generic point of each set of Table~\ref{t:typeII} are the same as those of a generic element of $A$, that is, they are given in the column
``decomposition of $V$''. In particular, this matrix is not diagonalizable. For $\so(4)$ on the hypersurface
$4qc_1=-(q^2-1)^2s_1^2$ the matrix $K_1$ has one eigenvalue with two blocks of size $2$, and $K_2$
is diagonalizable. At the intersection of the two hypersurfaces, $K_1$ and $K_2$ each have one eigenvalue with two
blocks of size $2$. On the sets of the first type the new solutions are $K$ and $K(s)$ for $\so(6)$ at $c_2=c_3$,
$K(s)$ for $\osp(2|4)$ at $c_2=c_3$ and $K(s,s_2)$ for $\so(5)$ under condition~\eqref{eq:sp}.

We formulate the results of Sections~\ref{ss:Ks} and~\ref{ss:special} as a theorem.

\begin{theorem}\label{th:cartan}
Let the Cartan term~\eqref{eq:Xs} be added to the generator of one of the nodes. For generic values of the
parameters, the dimensions of the centralizers and all solutions of the reflection equation are given in Table~\ref{t:resS} (for the rows with
the condition $c_2=c_3$, at a generic point of this hypersurface), and the
nonscalar solutions are written out in Section~\ref{ss:Ks}. If Cartan terms with independent parameters are added to the
generators of several nodes, then for generic values of the parameters the centralizer is nonscalar only for $\so(4)$. In addition, it is nonscalar for $\so(5)$ at a generic point of the hypersurface
$c_1=-(q-1)^2s^2$, and there the solutions of the reflection equation are $t\cdot1$ and $tK(s,s_2)$. There are no other special
values of the first type besides $c_2=c_3$ for $\so(6)$ and $\osp(2|4)$, including inside the sets of special values. All nonscalar solutions are invertible for all values of the parameters of the Cartan terms and nonzero $c_i$.
\end{theorem}

\begin{propn}\label{p:typeII}
In the cases of Table~\ref{t:dec}, for generic values of the parameters the centralizer $A$ is commutative and semisimple,
and $V$ decomposes over $\mathcal B$ without multiplicities, as indicated in Table~\ref{t:dec}. The sets of special values of the second
type in these cases (for $\so(5)$ with terms on both nodes, inside~\eqref{eq:sp}) are exhausted by
Table~\ref{t:typeII}: these are hypersurfaces and, for $\so(4)$ with terms on both nodes, their intersection. At a generic
point of each of these sets $A$ is commutative, the decomposition of $V$ is given in Table~\ref{t:typeII}, and the solutions
of the reflection equation in the centralizer are obtained by substituting the condition into the solutions at a generic point of the ambient
set. There are no other solutions. In the remaining cases (the centralizer at a generic point is scalar) there are no special values except
the sets of the first type.
\end{propn}

\subsection{Classical limit}\label{ss:cl}
The subalgebras $\kk$ and $\kk_\sigma$ are constructed as the graded closures of the generators under the supercommutator
(the even and odd parts of the generators are taken separately). The computations were carried out in general form:
$c_i$, $\sigma$ and $\sigma_2$ (the parameter of the term on $\alpha_2$ when the terms are on both nodes) remain symbols, and the ranks
and dimensions are found over the field of rational functions in them, that is, the results refer to generic values of the
parameters. We denote by $[\kk,\kk]$ the derived algebra, by $H_1(\kk)=\kk/[\kk,\kk]$ the first homology,
by $Z(\kk)$ the center of the algebra $\kk$, by $C(\cdot)$ the even centralizer in $\mathrm{End}(V)$, and by
$\operatorname{rad}C(\cdot)$ its radical (it is found by the lemma of Section~\ref{ss:special}). The generators of $\kk$ and
$\kk_\sigma$, the conditions on the parameters and the limits of the $K$-matrices are obtained from the quantum ones by the substitution $s=\sigma/(q-1)$
(and $s_i=\sigma_i/(q-1)$) and passing to the limit $q\to1$. The generator $T=q^{h}$ turns into $h$ here (the limit of
$(T-1)/(q-1)$). The results are collected in Tables~\ref{t:cl1}, \ref{t:cl2} and~\ref{t:cl3}.

\begin{table}[H]
\centering\small
\setlength{\tabcolsep}{4pt}
\begin{tabular}{@{}llcccccc@{}}
\toprule
algebra & Cartan term on & $\dim\g$ & $\dim\kk$ & $\dim\kk_\sigma$ & $\dim[\kk,\kk]$ & $\dim H_1(\kk)$ & $x_\alpha\in[\kk,\kk]$\\
\midrule
$\so(3)$ & $\alpha_1$ & 3 & 1 & 1 & 0 & 1 & no\\
$\so(4)$ & $\alpha_1$ and $\alpha_2$ & 6 & 2 & 2 & 0 & 2 & no\\
$\so(5)$ & $\alpha_1$ & 10 & 4 & 4 & 3 & 1 & no\\
$\so(5)$ & $\alpha_2$ & 10 & 4 & 10 & 3 & 1 & yes\\
$\so(5)$, $c_1=-\sigma^2$ & $\alpha_1$ and $\alpha_2$ & 10 & 4 & 4 & 2 & 2 & no\\
$\so(6)$, $c_2=c_3$ & $\alpha_1$ & 15 & 7 & 7 & 6 & 1 & no\\
$\so(6)$, $c_2=c_3$ & $\alpha_2$ & 15 & 7 & 15 & 6 & 1 & yes\\
$\so(6)$, $c_2\neq c_3$ & no term & 15 & 15 & 15 & 15 & 0 & yes\\
$\osp(1|2)$ & $\alpha_1$ (odd) & 5 & 2 & 5 & 1 & 1 & no\\
$\osp(1|4)$ & $\alpha_1$ & 14 & 6 & 6 & 5 & 1 & no\\
$\osp(1|4)$ & $\alpha_2$ (odd) & 14 & 6 & 14 & 5 & 1 & yes\\
$\osp(2|2)$ & $\alpha_1$ (odd) & 8 & 4 & 8 & 3 & 1 & yes\\
$\osp(2|4)$, $c_2=c_3$ & $\alpha_1$ & 19 & 9 & 9 & 7 & 2 & no\\
$\osp(2|4)$, $c_2=c_3$ & $\alpha_2$ (odd) & 19 & 9 & 19 & 7 & 2 & yes\\
$\osp(2|4)$, $c_2\neq c_3$ & $\alpha_1$ & 19 & 9 & 19 & 8 & 1 & yes\\
\bottomrule
\end{tabular}
\caption{Classical limit: dimensions of $\kk$, $\kk_\sigma$, the derived algebra and $H_1(\kk)$. In the row $\so(5)$, $c_1=-\sigma^2$, $\kk$ denotes the subalgebra with the term $\sigma h_1$ on $\alpha_1$,
$\kk_\sigma$ denotes the same subalgebra with the added term $\sigma_2h_2$ on $\alpha_2$, and $x_\alpha=x_2$.}\label{t:cl1}
\end{table}

\begin{table}[H]
\centering\small
\setlength{\tabcolsep}{3pt}
\begin{tabular}{@{}llccccc@{}}
\toprule
algebra & Cartan term on & $\dim Z(\kk)$ & $Z(\kk)\oplus[\kk,\kk]=\kk$ & $\dim C(\kk)$ & \hd{$\dim$\\ $C(\kk_\sigma)$} & \hd{limit\\ of $K$}\\
\midrule
$\so(3)$ & $\alpha_1$ & 1 & yes & 3 & 3 & yes\\
$\so(4)$ & $\alpha_1$ and $\alpha_2$ & 2 & yes & 4 & 4 & yes\\
$\so(5)$ & $\alpha_1$ & 1 & yes & 3 & 3 & yes\\
$\so(5)$ & $\alpha_2$ & 1 & yes & 3 & 1 & ---\\
$\so(5)$, $c_1=-\sigma^2$ & $\alpha_1$ and $\alpha_2$ & 1 & no & 3 & 3 & yes\\
$\so(6)$, $c_2=c_3$ & $\alpha_1$ & 1 & yes & 3 & 3 & yes\\
$\so(6)$, $c_2=c_3$ & $\alpha_2$ & 1 & yes & 3 & 1 & ---\\
$\so(6)$, $c_2\neq c_3$ & no term & 0 & yes & 1 & 1 & ---\\
$\osp(1|2)$ & $\alpha_1$ & 1 & no & 2 & 1 & ---\\
$\osp(1|4)$ & $\alpha_1$ & 1 & no & 3 & 3 & yes\\
$\osp(1|4)$ & $\alpha_2$ & 1 & no & 3 & 1 & ---\\
$\osp(2|2)$ & $\alpha_1$ & 1 & no & 2 & 1 & ---\\
$\osp(2|4)$, $c_2=c_3$ & $\alpha_1$ & 1 & no & 3 & 3 & yes\\
$\osp(2|4)$, $c_2=c_3$ & $\alpha_2$ & 1 & no & 3 & 1 & ---\\
$\osp(2|4)$, $c_2\neq c_3$ & $\alpha_1$ & 1 & yes & 2 & 1 & ---\\
\bottomrule
\end{tabular}
\caption{Center, even centralizers and the limit of $K(s)$. The last column indicates whether the limit of $K(s)$ as $q\to1$, $(q-1)s\to\sigma$
commutes with $\kk_\sigma$ (for $\so(4)$, the limits of $K_1$ and $K_2$). In the row $c_1=-\sigma^2$, this is
the limit of $K(s,s_2)$, $(q-1)s_2\to\sigma_2$.}\label{t:cl2}
\end{table}

For $\so(4)$ with a Cartan term only on $\alpha_1$ or only on $\alpha_2$, all the quantities of Tables~\ref{t:cl1} and~\ref{t:cl2}
are the same as in the row $\so(4)$ with terms on both nodes. The radical of $C(\kk_\sigma)$ is zero in all rows of Table~\ref{t:cl2}. The radical of $C(\kk)$ is nonzero in the rows $\osp(1|2)$ ($\dim\operatorname{rad}C(\kk)=1$), $\osp(1|4)$
(dimension $2$), $\osp(2|4)$ at $c_2=c_3$ (dimension $2$) and $\so(5)$, $c_1=-\sigma^2$ (dimension $2$). In the
remaining rows it is zero.

\textbf{Limit of the special values.} At a generic point of the limit of each set of Table~\ref{t:typeII} (the limit of the condition
is given in the same table), we found $\kk_\sigma$, $C(\kk_\sigma)$ and the decomposition of $V$ over $\kk_\sigma$. The results are in
Table~\ref{t:cl3}. The two rows of Table~\ref{t:typeII} marked ``none'' have no limit. The dimensions of
$\kk_\sigma$ and $C(\kk_\sigma)$ at a generic point of each limit are the same as at a generic point of the ambient set. Two
hypersurfaces of Table~\ref{t:typeII} with the same limit condition (for $\so(3)$, for $\so(5)$ with a term
on $\alpha_1$ and for $\so(6)$) turn into the same set, and at its generic point $\dim\operatorname{rad}C(\kk_\sigma)=2$,
while $V$ is indecomposable over $\kk_\sigma$. In the remaining rows the decomposition of $V$ over $\kk_\sigma$ is the same as over $\mathcal B$ in Table~\ref{t:typeII}.
The limits of the $K$-matrices at a generic point of each limit commute with $\kk_\sigma$, and the Jordan blocks of the limit of $K(s)$ (respectively
$K(s,s_2)$ or $K$, and for $\so(4)$ of the limit of $K_1K_2$) are the same as those of a generic element of $C(\kk_\sigma)$, that is, they are given
in the last column of Table~\ref{t:cl3}. In all rows of Table~\ref{t:cl3} except $\osp(2|2)$, the limit
condition is exactly the condition under which, for each node $\alpha$ whose parameters enter this condition,
the generator $x_\alpha+\sigma h_{\ta}$ has the form $e_\alpha-a^2f_\alpha-a[e_\alpha,f_\alpha]$ of the generators of the electrical
Lie algebra~\cite[Section~3.1]{AMS2} (see also~\cite{LP,BGG}), where $a=\pm\sigma$. For $\so(3)$ and the node $\alpha_2$ in $\so(5)$, where
$[[e_\alpha,f_\alpha],e_\alpha]=e_\alpha$, in this formula $f_\alpha$ is replaced by $2f_\alpha$, and $a=-\sigma/2$. In this case
$x_\alpha+\sigma h_{\ta}$ is nilpotent, whereas for generic values of the parameters it is semisimple.
For $\osp(2|4)$ without a Cartan term at $c_2=c_3$, the picture in the classical limit
is the same as in the quantum case: $\dim C(\kk)=3$, $\dim\operatorname{rad}C(\kk)=2$, and $V$ is indecomposable over $\kk$,
$(2|4)_{2,2,1,1}$.

\begin{table}[H]
\centering\small
\setlength{\tabcolsep}{4pt}
\begin{tabular}{@{}llccl@{}}
\toprule
algebra, Cartan term on & limit of the condition & $\dim C(\kk_\sigma)$ & $\dim\operatorname{rad}C(\kk_\sigma)$ & $V$ over $\kk_\sigma$\\
\midrule
$\so(3)$, $\alpha_1$ & $c_1=-\sigma^2/2$ & 3 & 2 & $3_{3}$\\
$\so(4)$, $\alpha_i$ or $\alpha_1,\alpha_2$ & $c_i=-\sigma_i^2$ & 4 & 2 & $2_{2}+2_{2}$\\
$\so(4)$, $\alpha_1$ and $\alpha_2$ & $c_1=-\sigma_1^2$, $c_2=-\sigma_2^2$ & 4 & 3 & $4_{3,1}$\\
$\so(5)$, $\alpha_1$ & $c_1=-\sigma^2$ & 3 & 2 & $5_{3,1,1}$\\
$\so(5)$, \eqref{eq:sp}, $\alpha_1$ and $\alpha_2$ & $c_2=-\sigma_2^2/2$ & 3 & 1 & $1+4_{2,2}$\\
$\so(6)$, $c_2=c_3$, $\alpha_1$ & $c_1=-\sigma^2$ & 3 & 2 & $6_{3,1,1,1}$\\
$\osp(1|4)$, $\alpha_1$ & $c_1=-\sigma^2$ & 3 & 1 & $(1|2)+(0|2)_{2}$\\
$\osp(2|2)$, no term & $c_1=c_2$ & 2 & 1 & $(2|2)_{2,1,1}$\\
$\osp(2|4)$, $c_2=c_3$, $\alpha_1$ & $c_1=-\sigma^2$ & 3 & 1 & $(2|2)+(0|2)_{2}$\\
\bottomrule
\end{tabular}
\caption{Classical limit of the special values of the second type: the limit of the condition from Table~\ref{t:typeII},
the centralizer $C(\kk_\sigma)$ and the decomposition of $V$ over $\kk_\sigma$ at a generic point of the limit. For $\so(3)$, $\so(5)$ with a term on $\alpha_1$ and
$\so(6)$, both hypersurfaces of Table~\ref{t:typeII} turn into this limit.}\label{t:cl3}
\end{table}

\subsection{Observation}\label{ss:hyp}
A Cartan term on a node $\alpha$ is called \emph{admissible} if there exists a nonscalar $K$-matrix with it.

\textbf{Observation.} In all cases of Table~\ref{t:cl1} the subalgebra $\kk_\sigma$ differs from $\g$ if and only
if $x_\alpha$ is even and $x_\alpha\notin[\kk,\kk]$, that is, the class of $x_\alpha$ in $H_1(\kk)=\kk/[\kk,\kk]$
is nonzero. Moreover, the Cartan term on $\alpha$ is admissible if and only if $\kk_\sigma\neq\g$
(Tables~\ref{t:resS} and~\ref{t:cl1}, and for the row $\so(5)$, $c_1=-\sigma^2$, Section~\ref{ss:special}).

\vspace{0.5cm}

\textbf{Acknowledgements.}

\medskip
This work was carried out at the Center of Pure Mathematics, MIPT.  \\
The authors E.~A.~Ryzhkin and V.~A.~Stukopin are financially supported by the Russian Science Foundation grant No.~26-11-00115.

\end{document}